\documentclass{svjour3}

\expandafter\let\csname equation*\endcsname\relax
\expandafter\let\csname endequation*\endcsname\relax

\pdfoutput=1

\usepackage[T1]{fontenc}
\usepackage{graphicx}
\usepackage{caption}
\usepackage{subcaption}
\usepackage{color}
\usepackage[backref=none]{hyperref}
\hypersetup{
  colorlinks,
  citecolor=red,
  filecolor=black,
  linkcolor=blue,
  urlcolor=magenta
}
\usepackage{amsmath}
\usepackage{amssymb}
\usepackage{lineno}
\usepackage{cite}
\usepackage{xparse}

\allowdisplaybreaks[1]

\DeclareMathSymbol{\Gamma}{\mathalpha}{operators}{0}
\DeclareMathAlphabet{\mathpzc}{OT1}{pzc}{m}{it}

\NewDocumentCommand\tauhit{m g}{%
  \IfValueF{#2}{\tau_{\textsc{#1}}}%
  \IfValueT{#2}{\tau_{\textsc{#1}}(#2)}%
}
\NewDocumentCommand\taumft{m g}{%
  \IfValueF{#2}{\overline{\tau}_{\textsc{#1}}}%
  \IfValueT{#2}{\overline{\tau}_{\textsc{#1}}(#2)}%
}
\NewDocumentCommand\fmft{m g}{%
  \IfValueF{#2}{\overline{f}_{\textsc{#1}}}%
  \IfValueT{#2}{\overline{f}_{\textsc{#1}}(#2)}%
}
\NewDocumentCommand\dG{g}{%
  \IfValueT{#1}{\partial \Gamma_{\textsc{#1}}}%
  \IfValueF{#1}{\partial \Gamma}%
}
\NewDocumentCommand\Gv{g}{%
  \IfValueT{#1}{\Gamma_{\textsc{#1}}}%
  \IfValueF{#1}{\Gamma}%
}
\NewDocumentCommand\Mf{g}{%
  \IfValueT{#1}{\mathpzc{M}_{\textsc{#1}}}%
  \IfValueF{#1}{\mathpzc{M}}%
}
\NewDocumentCommand\xf{g}{%
  \IfValueT{#1}{\mathpzc{#1}}%
  \IfValueF{#1}{\mathpzc{x}}%
}
\NewDocumentCommand\Mm{g}{%
  \IfValueT{#1}{\mathnormal{M}_{\textsc{#1}}}%
  \IfValueF{#1}{\mathnormal{M}}%
}
\NewDocumentCommand\xm{g}{%
  \IfValueT{#1}{\mathnormal{#1}}%
  \IfValueF{#1}{\mathnormal{x}}%
}
\NewDocumentCommand\num{g}{%
  \IfValueT{#1}{\nu_{\textsc{#1}}}%
  \IfValueF{#1}{\nu}%
}

\newcommand{\cnum}[2][bp]{\nu_{\textsc{#1}|#2}}
\newcommand{\expect}[2][\mu]{\mathbb{E}_{#1}(#2)}
\newcommand{\bb}[2]{\mathbb{#1}^{#2}}

\newcommand{\aref}[1]{appendix~\ref{#1}}
\newcommand{\sref}[1]{section~\ref{#1}}
\newcommand{\fref}[1]{figure~\ref{#1}}

\newcommand{\Pceq}[1][\beta]{P^{(\textsc{can})}_{#1}}
\newcommand{\eb}{\epsilon_{\textsc{b}} }
\newcommand{\ep}{\epsilon_{\textsc{p}} }
\newcommand{\teb}{\widetilde{\epsilon_{\textsc{b}} }}
\newcommand{\tep}{\widetilde{\epsilon_{\textsc{p}} }}

\newcommand{\thetah}{\,\Theta_\textsc{h}\!}
\newcommand{\canavg}[2]{\langle #1 \rangle^{(\textsc{can})}_{#2}}

\journalname{Journal of Statistical Physics}

\begin{document}

\title
{On the limiting Markov process of energy exchanges in a
  rarely interacting ball-piston gas}
\titlerunning{Stochastic limit of a rarely interacting ball-piston
  gas}
\dedication{To David Ruelle and Yasha Sinai on the occasion of their
  80th birthdays\footnote{We are deeply honored and privileged to
    dedicate this paper to David Ruelle and Yasha Sinai,
    great founding fathers of rigorous statistical physics. Their
    works were fundamental in creating the subject and have profoundly
    changed the way people think about it. In particular, beyond its
    wider and
    deeper impact, Ruelle's 1998 work \cite{Ruelle:1998response}, as
    well as his Brussels lecture \cite{Ruelle:1998Brussels}, were key
    to reinvigorating the general interest toward understanding
    Fourier's law of heat conduction and raising hopes for a
    satisfactory answer to this difficult problem. The theory of
    hyperbolic billiards, initiated by Sinai, helped engineer models
    which offer the most promising perspectives in this direction. It
    is our hope this paper testifies to their influence.}}
\author{P\'eter~B\'alint\dag\ddag
  \and
  Thomas~Gilbert\S
  \and
  P\'eter~N\'andori\P$\|$
  \and
  Domokos~Sz\'asz\dag
  \and
  Imre~P\'eter~T\'oth\dag\ddag
}
\authorrunning{P. B\'alint \emph{et al.}}
\institute{
  \dag\
  Institute of Mathematics, Budapest University of Technology and
  Economics, Egry J\'ozsef u. 1, H-1111 Budapest, Hungary
  \at
  \ddag\
  MTA-BME Stochastics Research group, Egry J\'ozsef u. 1, H-1111
  Budapest, Hungary
  \at
  \S\
  Center for Nonlinear Phenomena and Complex Systems,
  Universit\'e Libre  de Bruxelles, C.~P.~231, Campus Plaine, B-1050
  Brussels, Belgium
  \at
  \P\
  Department of Mathematics, University of Maryland, 4176 Campus Drive,
  College Park, MD 20742, USA
  \at
  $\|$\
  Courant Institute of Mathematical Sciences, New York University, New
  York, NY 10012 USA
}

\date{Version of \today}

\maketitle

\begin{abstract}
  We analyse the process of energy exchanges generated by the
  elastic collisions between a point-particle, confined to a
  two-dimensional cell with convex boundaries, and a `piston', i.e.~a
  line-segment, which moves back and forth along a one-dimensional
  interval partially intersecting the cell. This model can be
  considered as the elementary building block of a spatially extended
  high-dimensional billiard modeling heat transport in a class of
  hybrid materials exhibiting the kinetics of gases and
  spatial structure of solids. Using heuristic arguments and numerical
  analysis, we argue that, in a regime of rare interactions, the
  billiard process converges to a Markov jump process for the
  energy exchanges and obtain the expression of its generator.
\end{abstract}

\section{Introduction \label{sec:intro}}

Fourier's law of heat conduction \cite{Fourier:1822}, according to
which the heat current in a material is proportional to the gradient
of its local temperature, has over the last two centuries proved a powerful
phenomenological tool for describing the process of energy transfer in
physical systems. Yet, in spite of being well understood at a
macroscopic level, the derivation of this law from a microscopic point
of view arguably remains one of mathematical physics' great challenges.
Thus, in their millenium review, Bonetto \emph{et al.}
\cite{Bonetto:2000p13477}, after offering ``a selective overview of
the current state of our [then] knowledge (more precisely of our
ignorance) regarding the derivation of Fourier's Law'', proceeded to
the observation that
\begin{quote}
  ``There is however at present no rigorous mathematical derivation of
  Fourier's law [\dots] for any system (or model) with a
  deterministic, e.g. Hamiltonian, microscopic evolution.''
\end{quote}
An intriguing review of Fourier's work and its influence with a
detailed chronology of some of the main developments can be found in
reference \cite{Narasimhan:1999FourierHistory}. Other useful reviews
devoted to the non-equilibrium statistical mechanics of
low-dimensional systems include references \cite{Lepri:2003p8721} and
\cite{Dhar:2008p539}. 

Building upon the earlier work of Bunimovich \emph{et al.}
\cite{Bunimovich:1992p621}, Gaspard and Gilbert then set out in 2008
\cite{Gaspard:2008PRL101} to consider the \emph{regime of rare
  interactions} of a class of models, henceforth referred to as
GG-models, which, from
the point of view of ergodic theory \cite{Bunimovich:1992p621}, are
intermediate between the gas of hard balls and the periodic Lorentz
gas. In general, the GG-model can be thought of as a billiard chain (or
$\bb{Z}{d}$-network in dimension $2 \leq d \leq n$) of $n$-dimensional
cells with (semi-) dispersing walls, each containing a single ball particle
trapped inside it. Ball particles are moreover let to interact among
neighbours under the control of the geometry of the interface between
their respective cells. Hence by excluding mass transport the model
focuses on energy transport solely. 

A standard two-step strategy for analysing the process of heat
transport is to first identify the conditions under which a mesoscopic
description can be attained from the microscopic one, and second, by
taking the hydrodynamic limit of the mesoscopic process to
obtain, in the diffusive scaling, the heat equation at the
macroscopic scale. The completion of this step also implies gaining an
analytical form for the coefficient of heat conductivity.

The present work does not deal with the analysis of the second step of
this strategy. In particular, we will not address the precise form of
the coefficient of thermal conductivity associated with our model. We
note, however, that Sasada, inspired by Stefano Olla's remarks on the
results announced in reference \cite{Gaspard:2008PRL101},  recently
reported \cite{Sasada:2015} that the coefficient of heat conductivity
figuring in the papers of Gaspard and Gilbert, see in particular
reference \cite{Gaspard:200811P021}, corresponds to the contribution
to the heat conductivity from the static correlations alone, while the
true transport coefficient should also include a contribution from
dynamic correlations. Whereas the latter contribution appears to be very
small in comparison to the former, it does not vanish. This issue,
which deserves further clarification, will be the subject of future
work. Another important remark here is that for handling the second
step of the strategy for GG-models the necessary spectral gap has
already been obtained  in reference \cite{Sasada:2015Gap}; see also
reference \cite{Grigo:2012MixingRate}. The analogous lower estimate of
the spectral gap for the model to be introduced below is so far an
open question.

We mention here that progress in this area has been paralleled by
interesting prospects aiming at understanding the emergence of
Fourier's law when small noise is added to simple deterministic
models, e.g.~to weakly anharmonic crystals \cite{Liverani:2012Toward,
  Bricmont:2007Towards}, or harmonic oscillators
\cite{Bernardin:2005Fourier}. Such systems have also been considered
in regimes of rare interactions \cite{Keller:2009p2845,
  Huveneers:2012Energy}  Another line of research aims at 
studying energy exchanges between the cells of regular lattices
through the mediation of point particles
\cite{Eckmann:2006Nonequilibrium}; see also reference
\cite{Li:2016Local} for recent progress based on the seminal work of
Kipnis \emph{et al.} \cite{Kipnis:1982Heat}. A mixture of those two
approaches was treated in reference \cite{Bonetto:2004Fourier}. 

Going back to the first step of the strategy described above, in the
GG-model, the reduction from the deterministic dynamics at the
microscopic scale to a stochastic process at the mesoscopic scale,
emerges in the \emph{rare interaction limit} (RIL), out of a two-stage
relaxation process. In the RIL, the rate of interactions between
neighbouring ball particles is arbitrarily low, which implies that a
form of local equilibrium is achieved on the scale of every cell, each
kept at constant energy. The  convergence to local equilibrium is
indeed controlled by the rate of collisions between a particle and
the walls of its cell, which can be made much larger than the rate of
binary collisions, i.~e. of collisions between two neighbouring
particles. Since averaging takes place, energy transfers between
neighbouring particles behave stochastically, with every cell of the
system acting as a fundamental unit whose state is specified by the
energy of its ball particle. In effect, the RIL yields for the
Hamiltonian kinetic equation of any finite subsystem the generator of
a Markov jump process for the energies of the ball particles.

The GG-model can consist of networks of two-dimensional particles
(discs) confined to identical cells \cite{Gaspard:2008NJP3004}. It can
also consist of networks of three (or higher)-dimensional particles (spheres)
confined to identical cells \cite{Gaspard:2012p26117}. One can think
of extensions of such models with several particles trapped in every
cell; their numbers may be identical or vary from cell to cell and the
cells of the network may no be all identical. The dimensionality of
the dynamics would then vary from cell to cell. One main goal of the
present paper is to introduce yet another version of the  GG-model,
the simplest of its class, for which we deem the (mathematically
rigorous) completion of the first step of the GG-strategy realistic.

Let us explain why we think that such a modification is necessary.  The
`simplest' task in the original GG-programme is the treatment of a
two-cell system with two interacting discs. This is actually
isomorphic to a 4-dimensional semi-dispersing billiard. However,
statistical properties of higher-dimensional (larger than two)
billiards are so far understood exclusively for finite-horizon
strictly dispersing billiards \cite{Balint:2008p1309},
and even then only under the notorious `complexity hypothesis'.

It was therefore suggested in reference~\cite{PajorGyulai:2010p17699}
that by exploiting the RIL feature of the approach of
reference~\cite{Gaspard:2008PRL101}, a way out of this
high-dimensional quagmire would be to apply the method of `standard pairs'
\cite{Chernov:2009brownian}, a very efficient tool stemming from the
development of Markov approximation methods, which permits the use
of statistical properties of lower dimensional projections of the
model---in this case of well-understood planar Sinai-billiards. Be
that as it may, this idea led to further technical difficulties. 

For this reason we introduce below a \emph{ball-piston model}, which
belongs to the class of GG-models, but for which the dimension of the
(simplest) isomorphic semi-dispersing billiard shrinks from 4 to 3. In this
model, a disc particle caged in a two-dimensional cell with dispersing
walls is let to interact with a `piston' moving in a one-dimensional
interval. This is isomorphic to a 3-dimensional semi-dispersing
billiard. We believe that coping with its RIL is already a realistic
question, which is the subject of a longer project, involving four of
us, and still in progress;  as to the first related publication, see
\cite{balint:2015winp}. Suffice it to say here that the application of
the method of standard pairs to the ball-piston model is met with
difficulties similar to those alluded to above. However, in this model
their occurrence can be tracked and explicitly described and
ultimately, as we hope, treated. 

In this respect, it is worth mentioning that the method of standard
pairs  has already been successfully applied to models of heat
conduction. Thus, in reference \cite{Dolgopyat:2011Energy} a model of
weakly interacting Anosov flows was considered for which the
appropriate long time limit is a chain of interacting diffusion
processes rather than the kind of  Markov jump processes which are
expected to arise in rarely interacting systems. More recently, in
reference \cite{Dolgopyat:2016Encounter}, the authors obtained an
exponential limit law for the first encounter of two small discs on a
planar Sinai billiard table, a result somewhat closer to our program.

A further aim of the present work is the calculation of several
relevant characteristics of the ball-piston model, in particular:
\begin{itemize}
\item[(i)]
  The unconditional and conditional mean free times of binary
  collisions (conditioning on the outgoing energy partition between the
  two particles after a binary collision);
\item[(ii)]
 The transition kernel of the Markov jump process expected to emerge
 in the RIL;
\item[(iii)]
  The restriction of the invariant measure  to the binary
  collision surface at fixed energy.
\end{itemize}
Since some of these calculations are new, even as to their mathematical
content, we endeavour to formulate our arguments in a language which
we hope will be accessible to both mathematicians and
physicists.

Finally, we describe a computational test of the conjecture that the
emerging mesoscopic limit of the rarely interacting ball-piston model
is indeed the Markov jump process we claim it is.  This conjecture relies
on the assumption that, in the RIL, two successive binary collisions
are separated by enough wall collision events that averaging takes
place. We test this by considering several outgoing laws
from a binary collision and compare the ingoing laws at the next
binary collision event to the relevant equilibrium law by using the
Kullback-Leibler divergence \cite{Kullback:1951information}. We
provide numerical evidence that reducing the rate of binary collisions
yields a limiting distribution of energy exchanges consistent with the
expected result.

The paper is organised as follows. The ball-piston model is introduced
in \sref{sec:minmodel}. Section \ref{sec:mft} is devoted to the
calculation of the ball-piston mean free time and collision rate. In
\sref{sec:condmft} we introduce the notion of conditional mean free
time and proceed to its calculation. The derivation of the transition kernel of the
conjectured Markov jump process is provided in \sref{sec:stoch},
together with a description of our numerical test of the validity of
the conjecture, as well as a discussion of numerical
results. Concluding remarks are given in \sref{sec:conc}. Specific
calculations are provided in \aref{app:vol} (volume integrals),
\aref{app:wall} (wall collision frequencies) and \aref{app:cnuep}
(restriction of the invariant measure  to the binary collision surface
at fixed energy).

\section{Minimal ball-piston model \label{sec:minmodel}}

\begin{figure}[htb]
  \centering
  \includegraphics[width=0.9\textwidth]
  {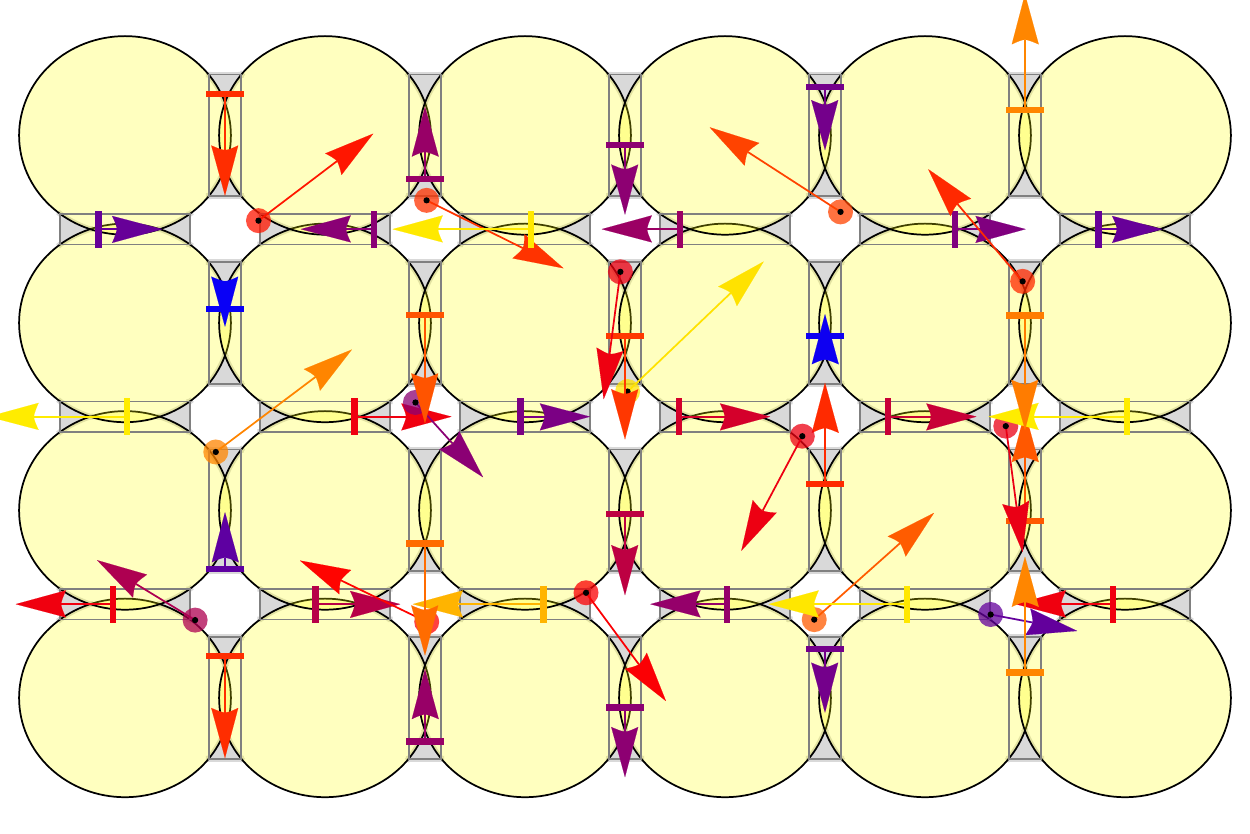}
  \caption{Ball-piston gas in a random configuration of the positions
    and velocities of the balls and pistons. The arrows' lengths and colours
    reflect the magnitude of the corresponding particle or piston's
    kinetic energy (blends of blue for low energy values, red and
    yellow for high energy values). Here periodic boundary
    conditions are applied: the piston cells of the right-most column
    are identical to those of the left-most column, and similarly for
    the upper and lower rows.}
  \label{fig:bpmanycells}
\end{figure}

The ball-piston gas, shown in \fref{fig:bpmanycells}, is a
collection of alternating balls and pistons, arranged in a periodic
structure, with every particle confined to its own cell. Balls and
pistons are particles of two different types. On the one hand,
balls are point-particles with two degrees of freedom. They move in
two-dimensional closed cells whose boundaries are defined by
impenetrable circular obstacles placed at the vertices of a square
lattice. Pistons, on the other hand, have only one degree of
freedom. They are one-dimensional vertical or horizontal line-segments
that are allowed to move back and forth along perpendicular intervals
placed between two neighbouring ball cells. Whereas pistons are
unaffected by the presence of the circular obstacles in the ball
cells, they do interact elastically with balls whenever collisions
occur, thereby exchanging their horizontal or vertical velocity
components (both balls and pistons have unit masses). By choosing the
lengths of the piston cells large enough that  their extremities
lay inside the ball cells (one symmetrically on each side), we
allow for energy exchanges between every ball and piston pair, the
likelihood of which depends on the piston's penetration length into
the ball cell. In turn, whereas mass transport of either species is
prohibited by the confining walls in every cell, energy exchanges
between balls and pistons induce heat transport on the scale of the
ball-piston gas.

\begin{figure}[htb]
  \centering
  \begin{subfigure}[t]{0.61\textwidth}
    \includegraphics[width=\textwidth]
    {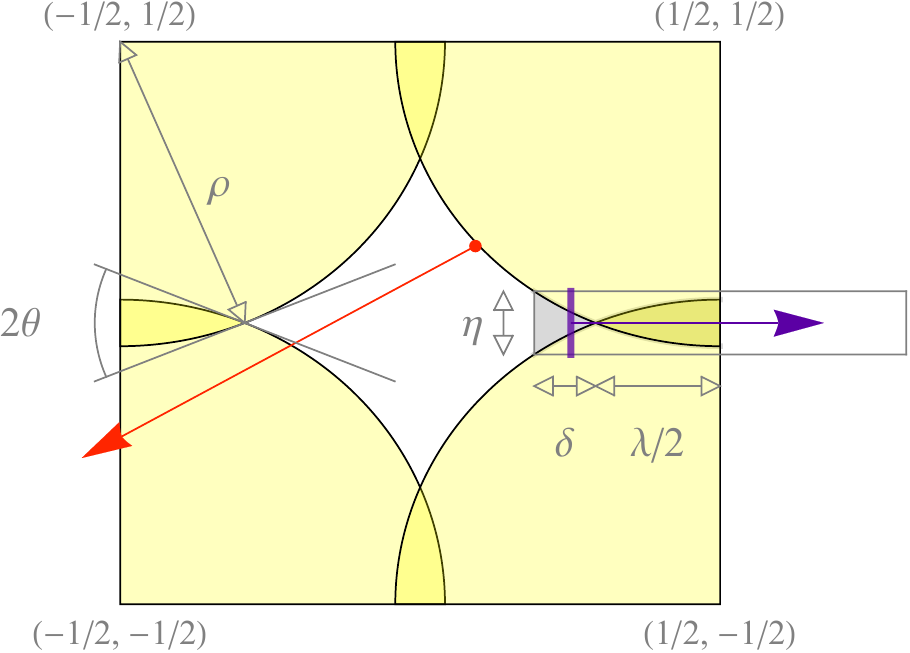}
    \caption{Minimal ball-piston model}
    \label{fig:modelsglcell}
  \end{subfigure}
  \begin{subfigure}[t]{0.37\textwidth}
    \includegraphics[width=\textwidth]
    {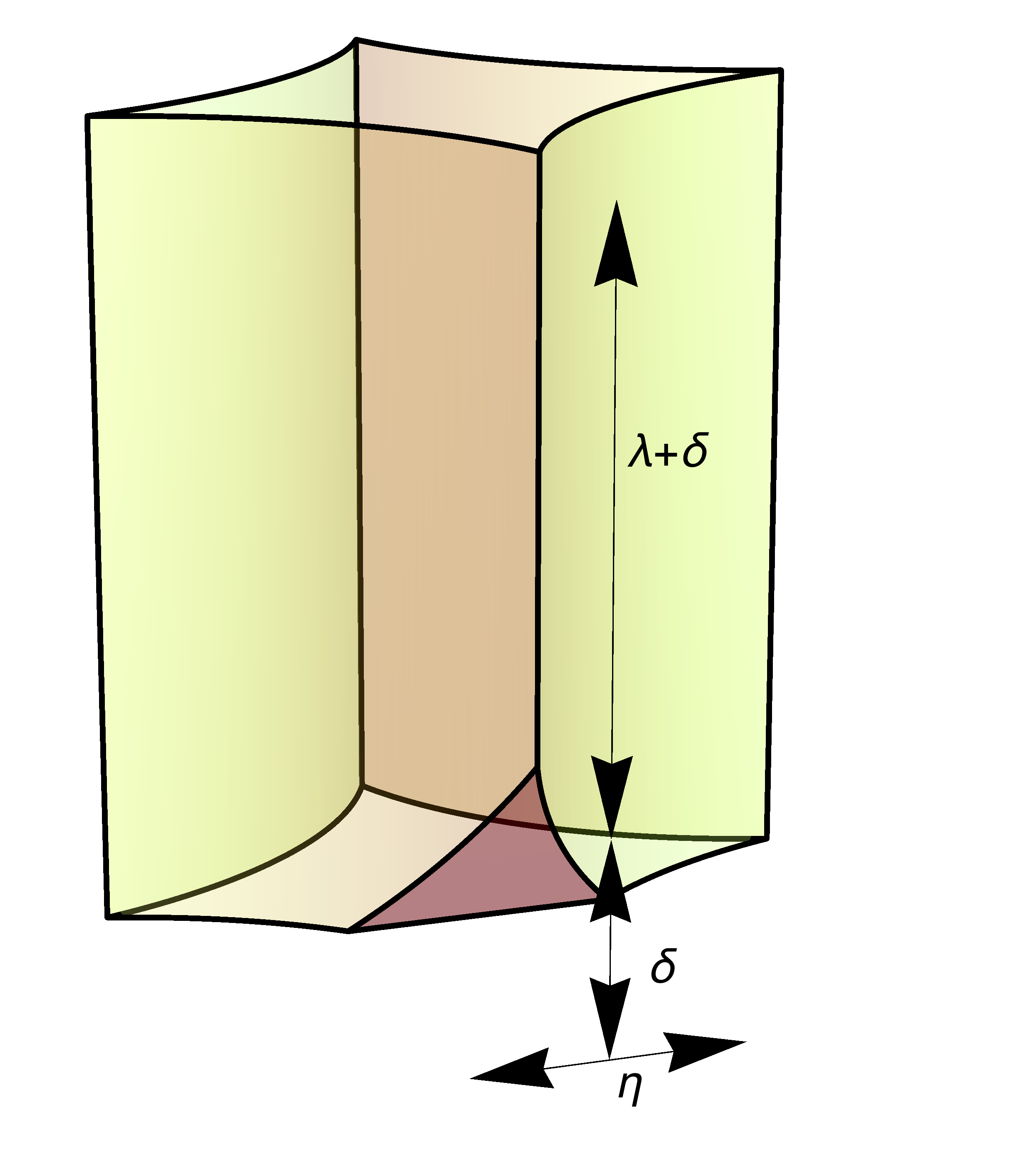}
    \caption{Three dimensional rendition}
    \label{fig:model3d}
  \end{subfigure}
  \caption{(a) The minimal ball-piston model consists of a single
    ball-piston pair, here showed in a random configuration. The
    relevant parameters are defined in the text. (b) Three-dimensional
    rendition of the billiard boundary.}
  \label{fig:model}
\end{figure}

The present study focuses on a minimal version of the ball-piston gas,
such as shown in \fref{fig:modelsglcell}. Here a single pair of
ball and piston---in this case a horizontally moving vertical
line-segment---is let to interact. This model can in fact be viewed as a
three-dimensional billiard, rendered in \fref{fig:model3d}: it is
equivalent to the free motion of a point-particle in a
three-dimensional cavity, undergoing elastic reflections upon its
boundary. The corresponding collision map is a four-dimensional
symplectic map.

The parameters relevant to the definition of the model are displayed
in \fref{fig:model}.  A point-particle (ball) of unit mass moves
freely in the interior of a cell (the ball cell) whose boundaries are
delimited by four arc-circles of common radius $\rho$, $1/2 < \rho <
1/\sqrt2$, centered at the four corners of a unit cell, and performs
elastic collisions with them. A vertical line segment of height $\eta$
and unit mass, which we call piston, moves horizontally back and forth
between the two edges of an interval of length $\lambda + 2 \delta$
(the piston cell), centered at the middle point of the cell's right
edge, where $\lambda = \sqrt{4\rho^2 - 1}$ measures the length  of the
interval between the two intersections of opposite discs (and such
that $\tan \theta = \lambda$). The parameter\footnote{The upper bound
  on $\delta$ is imposed so as to prevent overlap between the piston
  and a similar hypothetical vertical piston centered on either of the
  top and bottom edges of the cell, such as in the cells depicted in
\fref{fig:bpmanycells}.} $\delta$, $0<\delta < \rho/\sqrt{2} - \lambda/2$,
measures the length of penetration of the piston inside the ball cell and
therefore determines the possibility of interactions between the ball
and the piston. The height of the piston, $\eta = 1 - 2\sqrt{ \rho^2
  - (\lambda/2 +   \delta)^2}$, is such that, at its left-most
position, it lies inside the boundary of the ball cell. The positions
of the ball and piston must be initially chosen so that the ball is
located to the left of the piston; the ball cannot move passed the
piston.

We let $\Gv \subset \bb{R}{3}$ denote the three-dimensional
ball-piston configuration space and $\dG = \dG{bp} \cup \dG{bw} \cup
\dG{pw}$ its boundary, where $\dG{bp}$ is the surface of ball-piston
collisions, $\dG{bw}$ the surface of ball-wall collisions, and
$\dG{pw}$ the surface of piston-wall collisions. In
\fref{fig:model3d}, the first term refers to the slanted darker
surface of triangular shape, the second to the vertical walls, and the
third to the flat top and bottom walls. The phase space of the
billiard flow $\Mf$ is the product of $\Gv$ and $\bb{S}{2}$, the
sphere of unit radius in $\bb{R}{3}$, $\Mf = \Gv \times
\bb{S}{2}$.

A point $\mathbf{q} =
(q_1, q_2, q_3) \in \Gv$ specifies the horizontal and vertical
coordinates of the ball, $(q_1, q_2)$, and the piston's position,
$q_3$. They are such that:

\begin{equation}
  \label{eq:x1x2x3}
  \begin{split}
    &
    (q_1 \pm \tfrac{1}{2})^2 + (q_2 \pm \tfrac{1}{2})^2 \geq \rho^2\,,
    \\
    &
    \tfrac{1}{2}(1 - \lambda) -
    \delta \leq q_3 \leq \tfrac{1}{2}(1 + \lambda) + \delta\,,
    \\
    & q_1 \leq q_3\,.
  \end{split}
\end{equation}

\noindent The associated velocity vector is $\mathbf{v} = (v_1, v_2,
v_3) \in \bb{S}{2}$. The system's total (kinetic) energy is the sum
$\eb + \ep = \tfrac{1}{2}$ of the ball and piston energies, $\eb =
\tfrac{1}{2}(v_1^2 + v_2^2)$ and $\ep = \tfrac{1}{2}v_3^2$
respectively. The corresponding phase-space point is
denoted $\xf = \{\mathbf{q}, \mathbf{v}\} \in \Mf$.

 The phase space of the billiard map is denoted by $\Mm$
and is given by the product of $\dG$ and the set of vectors
$\mathbf{v}\in\bb{S}{2}$ whose scalar product with the unit
vector normal to the billiard surface at point $\mathbf{q}\in\dG$ is
non-negative. We have the decomposition $\Mm = \Mm{bp} \cup \Mm{bw}
\cup \Mm{pw}$. We write $\xm = \{\mathbf{q}, \mathbf{v}\} \in \Mm$ a phase
point of the billiard map.

The natural invariant measure of the flow is denoted by $\mu$ and
is normalised so that $\mu(\Mf) = 1$. Likewise the invariant measure
of the billiard map is denoted by $\num$ and such that $\num(\Mm) = 1$.

\section{Mean free times, collision frequencies and
  rates \label{sec:mft}} 

Let $S^t\xf$ denote the flow generated by the billiard dynamics on $\Mf$.
The first hitting time is the function

\begin{equation}
  \label{eq:hittime}
  \tauhit{} :\, \Mf\mapsto \bb{R}{+}:\, \xf \mapsto
  \inf\{t>0\,|\, S^t \xf \in \Mm\}
  \,.
\end{equation}

\noindent For $\xm\in\Mm$, $\tauhit{}(\xm)$ is the return time to the
billiard surface \cite{Cornfeld:book1982}, or free (flight) time
\cite{Chernov:1997p1}. Similarly, we define the ball-piston
free flight time to be the time separating two successive collisions
between the ball and piston,

\begin{equation}
  \label{eq:bpfreetime}
  \tauhit{bp}:\, \Mm{bp}\mapsto \bb{R}{+}:\, \xm \mapsto
  \inf\{t>0\,|\, S^t \xm \in \Mm{bp}\}
  \,.
\end{equation}

By ergodicity, the mean free time, which is defined to be the infinite
$n$ limit of the time to the $n$th collision (with any of the surface
elements of the three-dimensional billiard cavity) divided by the
number of collisions $n$, almost surely exists and is independent of
the initial condition if the latter is sampled with respect to the
natural invariant measure of the billiard flow. It is then equal to

\begin{equation}
  \label{eq:taumft}
  \taumft{} = \expect[\num]{\tauhit{}} \equiv \int_{\Mm} \tauhit{} (\xm)
  \,\mathrm{d}\num(\xm)
  \,,
\end{equation}

\noindent measured in terms of the natural invariant measure of the
billiard map on $\Mm$.

The ball-piston mean free time, i.e.~the average time separating
successive collisions of the ball-piston pair, is defined similarly to be

\begin{equation}
  \taumft{bp}
  = \expect[\num{bp}]{\tauhit{bp}}
    \equiv
    \int_{\Mm{bp}} \tauhit{bp}(\xm)
    \,\mathrm{d}\num{bp}(\xm),
    \label{eq:taubpmftint}
\end{equation}

\noindent where the measure $\num{bp}$ is $\num$ conditioned on
$\Mm{bp}$, $\num{bp}=\num{}(\Mm{bp})^{-1}\num{} |_{\Mm{bp}}$, which is
the natural invariant measure of the first return map from $\Mm{bp}$ to
itself.

As explained in references~\cite{Chernov:1991p204, Chernov:1997p1},
the presence of the hitting time under the integral in equations
\eqref{eq:taumft} and \eqref{eq:taubpmftint} has the effect
of lifting the integral on $\Mm$ to a measure on $\Mf$. We thus obtain an
explicit formula for the ball-piston mean free time by taking the
ratio between the normalising factors of the two invariant measures,
that of the flow to that of  the map restricted to $\Mm{bp}$. For the
billiard flow, we write $\mathrm{d}\mu(\xf) = c_{\mu} \,
\mathrm{d}\mathbf{q} \,  \mathrm{d}\mathbf{v}$, and, for the
conditional measure of the billiard map $\mathrm{d}\num{bp}(\xm) =
c_{\num{bp}} \mathrm{d}\mathbf{q} \, \mathrm{d}\mathbf{v} \,
(\mathbf{v}\cdot\mathbf{n})$, where $\mathbf{n}$ is the
($\mathbf{q}$-independent) unit vector normal to $\dG{bp}$,

\begin{equation}
  \label{eq:n}
  \mathbf{n} = \frac{1}{\sqrt2}(-1,0,1) \,.
\end{equation}

\noindent These normalising factors are, respectively,

\begin{align}
  \label{eq:cmu}
  c_{\mu}^{-1} &= \int_{\Gv} \mathrm{d}\mathbf{q}\,
  \int_{\bb{S}{2}} \mathrm{d}\mathbf{v}\, = 4 \pi |\Gv|
                 \,,
  \\
  \label{eq:cnu}
  c_{\num{bp}}^{-1} &= \int_{\dG{bp}} \mathrm{d}\mathbf{q}\,
  \int_{\bb{S}{2}: \, \mathbf{v}\cdot\mathbf{n} > 0}
                 \mathrm{d}\mathbf{v}\, (\mathbf{v}\cdot\mathbf{n}) =
                 \pi |\dG{bp}| \,,
\end{align}

\noindent where we have substituted $|\bb{S}{2}| = 4\pi$ and $|\bb{B}{2}| =
\pi$, the volume of the unit ball in $\bb{R}{2}$.

The ratio between equations \eqref{eq:cmu} and \eqref{eq:cnu} yields
the ball-piston mean free time, or mean return time to the ball-piston
collision surface,

\begin{equation}
  \label{eq:taubpmft}
  \taumft{bp}= \frac{c_{\num{bp}}}{c_{\mu}} = \frac{4 |\Gv|}{|\dG{bp}|}\,.
\end{equation}

\noindent This formula is but a special case of the well-known formula
for the mean-free time of three-dimensional billiards
\cite{Chernov:1997p1}.

The simple geometry of the minimal ball-piston model allows for an
explicit computation of the volume and surface integrals in equation
\eqref{eq:taubpmft}:

\begin{align}
  |\Gv|
  &=
    (\lambda + 2 \delta)
    \Big\{1 - \lambda - \rho^2
    \big[ \pi - 4 \arctan( \lambda) \big]\Big\}
    \cr
  &\quad
    - \frac{1}{24} \Big\{
    2 \delta ( \lambda + 4 \delta ) +
    \big[ 1 - \sqrt{ 1 - 4 \delta ( \lambda + \delta ) } \big]
    \big[ 2 + 4 \delta ( \lambda + \delta ) + 3 \lambda^2 \big]
    \Big\}
    \cr
  &
    \quad
    -
    \frac{1}{2} \rho^2 ( \lambda + 2\delta )
    \Big[ \arctan \lambda
    - \arctan \frac {\lambda + 2 \delta}
    {\sqrt{  1 - 4 \delta ( \lambda + \delta )}}
    \Big]\,,
    \label{eq:Gv}
  \\
  |\dG{bp}|
  &=
    \frac{1}{2\sqrt 2} \Big\{
    ( \lambda + 2 \delta )
    \big[ 2 - \sqrt{1 - 4 \delta  ( \delta + \lambda)} \big]
    - \lambda
    \Big\}
    \cr
  &\quad + \sqrt2 \rho^2 \Big[
    \arctan \lambda -
    \arctan \frac{\lambda + 2 \delta}
    {\sqrt{1 - 4 \delta  ( \delta + \lambda)}}
    \Big]
    \,;
    \label{eq:dGc}
\end{align}

\noindent see \aref{app:vol} for details.

The small $\delta$ regime, when ball-piston collisions become
arbitrarily rare,  is of particular interest, as discussed in
\sref{sec:stoch}. On the one hand, $\lim_{\delta\to0}|\Gv|$ is simply
the area \eqref{eq:areax1x2} multiplied by $\lambda$, the length of
the piston's interval of motion \eqref{eq:x1x2x3} in the limit
$\delta\to0$. On the other hand, the piston's height in this regime is
$\eta \simeq 2 \lambda \delta$, so that the region the piston can
penetrate inside the ball cell forms approximately a triangle of area
$\lambda \delta^2$ in the $(q_1, q_2)$ plane. As discussed in
\aref{app:vol}, this triangle is essentially the projection on the
ball cell of the Poincar\'e section of ball-piston
collisions. Accordingly, $|\dG{bp}|
\simeq \sqrt 2 \lambda \delta^2$. To leading order, the inverse of the
mean free time is therefore proportional to the parameter
squared\footnote{For $\rho = \tfrac{1}{2}$, however, we have $\lambda
  = 0$ so that, in the small $\delta$ regime,  $|\Gv| \simeq 2 \delta
  (1 - \tfrac{\pi}{4})$ and $|\dG{bp}| \simeq \tfrac{2\sqrt 2}{3}
  \delta^3$. The two limits $\rho\to\tfrac{1}{2}$ and $\delta\to0$ are
  therefore not interchangeable:
  \[
  \label{eq:limitsrhodelta}
  \lim_{\rho\to{1}/{2}} \lim_{\delta\to0} \delta^2 \frac{|\Gv|}{|\dG{bp}|} =
  \frac{1}{4 \sqrt{2}} (4 - \pi) \neq
  \lim_{\delta\to0} \lim_{\rho\to{1}/{2}} \delta^2 \frac{|\Gv|}{|\dG{bp}|} =
  \frac{3}{ 4 \sqrt {2}} ( 4 - \pi )\,.
  \]
},

\begin{equation}
  \label{eq:taubpmftdelta}
  \lim_{\delta\to0} (\taumft{bp} \delta^2)^{-1}
  = \frac{1}{2\sqrt2}
  \Big[
  1 - \lambda - \rho^2 \big( \pi - 4 \arctan \lambda \big)
  \Big]^{-1}
  \,,
\end{equation}

\noindent which, up to the prefactor, is the inverse of area
\eqref{eq:areax1x2}.

Since the number of collisions up to some time $t$ almost
surely increases like $t/\taumft{}$ as $t\to\infty$, it is natural
to call $\fmft{} \equiv \taumft{}^{-1}$ the collision frequency. However,
$\fmft{}$ can also be identified with a probability rate. Indeed,
ergodicity implies that $\fmft{}$ is also, for any time interval, the
expected number of collisions measured in that time interval
divided by its length. Thus, in particular, $\fmft{}$ can be expressed using
the probability to observe at least one collision up to time $t$, which is
$\mu \big(\{\xf \in \Mf \,|\, S^t \xf \cap \Mm \neq \emptyset\}\big)$. Namely,

\begin{equation}
  \label{eq:fmftrate}
  \fmft{} = \lim_{t\to0}
  \frac{1}{t}
  \mu \big(\{\xf \in \Mf \,|\, S^t \xf \cap \Mm \neq \emptyset\}\big)\,.
\end{equation}

By the same token, $\fmft{bp} \equiv\taumft{bp}^{-1}$ is the
ball-piston collision frequency, and also defines a probability rate
in the sense that

\begin{equation}
  \label{eq:fmftbprate}
  \fmft{bp} = \lim_{t\to 0}\frac{1}{t}
  \mu \big( \{\xf \in \Mf \,|\, S^t \xf \cap \Mm{bp} \neq \emptyset\} \big)\,.
\end{equation}

We emphasise that equations \eqref{eq:fmftrate} and
\eqref{eq:fmftbprate} justify referring to $\fmft{}$ and $\fmft{bp}$
as collision (probability) rates, regardless of the actual
distributions of the waiting times $\tauhit{}$ and $\tauhit{bp}$. In
particular, these quantities may have distributions which are not
exponential, and, accordingly, the collision event process may not be
Poisson.

Similar considerations apply to the ball-wall and piston-wall
collision events. We refer to \aref{app:wall} for a computation of the
corresponding return times.

\section{Conditional mean free times, collision frequencies and
  rates \label{sec:condmft}} 

Say the billiard flow is in a stationary regime and we are observing
the successive collision events that result in energy exchanges between
the ball and piston. Let us assume we are only interested in a
marginal set of events such that the ball-piston energy partition has
the value $\{\eb, \ep\}$, with $\eb + \ep = \tfrac{1}{2}$. Specifically, we ask:
for points $\xm\in\Mm$ whose velocity vectors $\mathbf{v} = (v_1, v_2,
v_3)$ are such that $\tfrac{1}{2}(v_1^2 + v_2^2) = \eb$ and
$\tfrac{1}{2}v_3^2 = \ep$, what is the corresponding mean free time? A
formula similar to equation \eqref{eq:taubpmftint} is obtained for
this quantity, which we
call the conditional mean free time:

\begin{equation}
  \label{eq:condbpmeanfreetime}
  \taumft{bp} (\ep) = \expect[\cnum{\ep}]{\tauhit{bp}}
  \equiv \int_{\Mm{bp}(\ep)}
  \tauhit{bp}(x) \, \mathrm{d}\cnum{\ep}(x)\,.
\end{equation}

\noindent where the measure $\cnum{\ep}$ is the measure $\num{bp}$
conditioned on the subset $\Mm{bp}(\ep) \subset \Mm{bp}$ of
phase-space points $\xm = \{\mathbf{q}, \mathbf{v}\}$ such that
$\mathbf{q}\in\dG{bp}$ and
$\mathbf{v} = (v_1,v_2,v_3) \in\bb{S}{2}$ with $v_3 = \pm \sqrt{2\ep}$
and $v_3 > v_1$ (consistent with $\mathbf{v}\cdot\mathbf{n}>0$).

To obtain an explicit formula, it is enough to transpose the
normalising factors in equations \eqref{eq:cmu} and \eqref{eq:cnu} to
the associated subvolumes of phase-space,

\begin{equation}
  \label{eq:taubpmftep}
  \taumft{bp}(\ep) = \frac{c_{\cnum{\ep}}}{c_{\mu_{\ep}}}\,.
\end{equation}

\noindent These normalising factors are computed as follows.
Note that the three-dimensional velocity vector $\mathbf{v}$ is
allowed to take values on two circles of radii $\sqrt{1 -
  2\ep}$ parallel to the plane $(v_1, v_2)$, whose two heights
correspond to the two allowed signs of the piston velocity $v_3$,

\begin{equation}
  \label{eq:v1v2v3ep}
  \begin{split}
    v_1 & = \sqrt{1 - 2 \ep} \cos \alpha\,,\\
    v_2 & = \sqrt{1 - 2 \ep} \sin \alpha\,,\\
    v_3 & = \sigma \sqrt{2\ep}\,,
  \end{split}
\end{equation}

\noindent where $\sigma=\pm1$; see \fref{fig:collsurface} in
\aref{app:cnuep}. The inverse of the factor $c_{\mu_{\ep}}$ is thus
the product of the volume $|\Gv|$ and twice the perimeter of the unit
circle,

\begin{equation}
  \label{eq:cmuep}
    c_{\mu_{\ep}}^{-1} = \int_{\Gv} \mathrm{d}\mathbf{q}\,
  \int_{2\times\bb{S}{1}} \mathrm{d}\mathbf{v}\, = 4 \pi |\Gv| \,,
\end{equation}

\noindent which turns out to be identical to $c_{\mu}$, equation
\eqref{eq:cmu}. This quantity must indeed be independent of the
parameter $\ep$, since only the orientations of the velocities are
relevant. This is, however, not so of the factor $c_{\cnum{\ep}}$,
which, as described in \aref{app:cnuep}, is found to be:

\begin{align}
  c_{\cnum{\ep}}^{-1}
  &=
    \int_{\dG{bp}} \mathrm{d}\mathbf{q}\,
    \int_{2\times \bb{S}{1}: \, \mathbf{v}\cdot\mathbf{n} > 0}
    \mathrm{d}\mathbf{v}\, (\mathbf{v}\cdot\mathbf{n}) \,,
    \cr
  &=
    4 \pi |\dG{bp}|
    \begin{cases}
      \tfrac{1}{\pi}
      \left[\sqrt{\tfrac{1}{2} - 2\ep} + \sqrt{\ep} \arcsin
        \sqrt{\frac{\displaystyle \ep}{\displaystyle \tfrac{1}{2} - \ep}}
      \right]  \,,
      & \ep < \tfrac{1}{4}\,,\\
      \tfrac{1}{2} \sqrt{\ep}\,,
      & \ep \geq \tfrac{1}{4}\,.
    \end{cases}
  \label{eq:cnuep}
\end{align}

To perform an actual measurement of the conditional mean free time
\eqref{eq:taubpmftep}, one must sample initial conditions with respect
to the density

\begin{equation}
  \label{eq:nuepdensity}
  c_{\cnum{\ep}} (\mathbf{v} \cdot \mathbf{n})_+
\end{equation}

\noindent on $\Mm{bp}$, where $(x)_+ = x$ if $x>0$, and $0$
otherwise; see \aref{app:cnuep} for details. It is, however, worth
noting the conditional mean free time \eqref{eq:taubpmftep} may be
computed most simply as follows. Consider an equilibrium time series
of the billiard dynamics and select the subset of ball-piston collision
events with energy partitions close to the desired one. Indeed,
definition \eqref{eq:condbpmeanfreetime} can immediately be extended
to arbitrary  energy intervals. Considering, in particular, the
interval $(\ep - \tfrac{1}{2} \varepsilon , \ep +
\tfrac{1}{2}\varepsilon)$ for small $\varepsilon>0$, we have the identity

\begin{equation}
  \label{eq:taumftepdelta}
  \taumft{bp}(\ep) = \lim_{\varepsilon\to0} \frac{\taumft{bp}(\ep -\tfrac{1}{2}
    \varepsilon ,\ep + \tfrac{1}{2}
    \varepsilon)}{\varepsilon}
  \,,
\end{equation}

\noindent which guarantees the convergence, as one decreases the
parameter $\varepsilon$, to the desired result of a measurement
performed on a coarser set.

The reason that makes the conditional mean free time
\eqref{eq:taubpmftep} particularly interesting is that its inverse
$\fmft{bp}{\ep}  \equiv \taumft{bp}(\ep)^{-1}$ can again be
viewed as a rate,

\begin{equation}
  \label{eq:fmftbpeprate}
  \fmft{bp}{\ep}  =
  \lim_{t\to 0}\frac{1}{t}
  \mu_{\ep} \big( \{\xf \in \Mf \,|\, S^t \xf \cap \Mm{bp} \neq \emptyset\} \big)\,.
\end{equation}

\noindent In the rare interaction limit $\delta\to 0$, we expect
the conditional distribution of $\tauhit{bp}(\ep)$ to indeed become
exponential, with rate $\fmft{bp}{\ep}$, consistent with the
expectation that the energy process converges to a Markov jump
process; see \sref{sec:stoch}.

Note that, by definition of $c_{\cnum{\ep}}$, we recover
the inverse of $c_{\num{bp}}$ after integrating the inverse of the
former quantity over the values of the
piston energy $\ep$, weighted by the density of
its marginal equilibrium distribution, a Beta distribution of shape
parameters $\tfrac{1}{2}$ and $1$\footnote{The integral of the surface element
  on $\bb{S}{2}$ along the two horizontal circles at heights
  $\pm\sqrt{2\ep}$ is $1/\sqrt{2\ep}$, which is the density of the
  Beta distribution of shape parameters $\tfrac{1}{2}$ and $1$, properly
  normalised.}. This relation implies a similar one between the
ball-piston conditional collision rate \eqref{eq:fmftbpeprate} and the
collision rate \eqref{eq:fmftbprate},

\begin{equation}
  \label{eq:intcnuep}
  \int_{0}^{\tfrac{1}{2}} \mathrm{d}\ep\,
  \frac{1}{\sqrt{2\ep}}
  \fmft{bp}{\ep} = \fmft{bp}\,.
\end{equation}

\noindent This does not contradict the identities $\fmft{bp}{\ep} =
\taumft{bp}(\ep)^{-1}$ and $\fmft{bp} = \taumft{bp}^{-1}$.

The product of the mean free time \eqref{eq:taubpmft} by the collision
rate \eqref{eq:fmftbpeprate} allows to define a dimensionless
energy-dependent collision frequency which is independent of the
billiard geometry,

\begin{equation}
  \label{eq:numft}
  \phi_{\textsc{bp}}(\ep) \equiv \taumft{bp} \fmft{bp}{\ep} =
  \begin{cases}
    \frac{4} {\pi}
    \left[
      \sqrt{\tfrac{1}{2} - 2\ep} + \sqrt{\ep} \arcsin
      \sqrt{\frac{\displaystyle \ep}{\displaystyle \tfrac{1}{2} - \ep}}
    \right]  \,,
    & 0 < \ep \leq \tfrac{1}{4}\,,\\
    2 \sqrt{\ep}\,,
    & \tfrac{1}{4} < \ep \leq \tfrac{1}{2}\,.
  \end{cases}
\end{equation}

\begin{figure}[htb]
  \centering
  \includegraphics[width=0.8\textwidth]
  {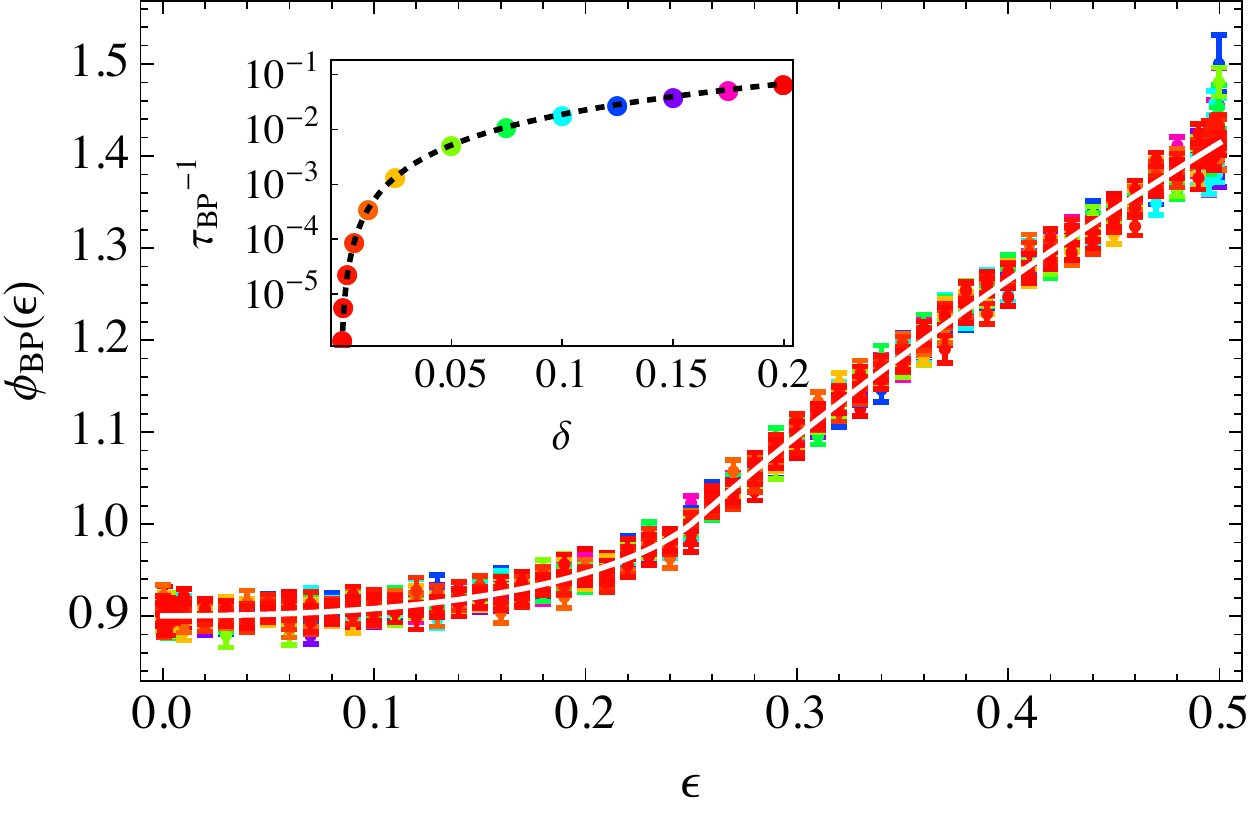}
  \caption{Numerical computations of the product of the mean free time
    by the conditional collision rate, compared with the rescaled
    ball-piston collision frequency, $\phi_{\textsc{bp}}(\ep)$, given
    by equation \eqref{eq:numft} (solid white curve), plotted as
    function of the piston energy $\ep = \epsilon$ (with the ball
    energy $\eb = \tfrac{1}{2} - \epsilon$). The parameter $\delta$
    takes on a number of different values which can be read off the
    horizontal axis of the graph in the inset. There we show the
    numerical computations of the inverse of the mean free time as
    function of $\delta$ in comparison to equation \eqref{eq:taubpmft}
    (dotted black curve). The colours of the data in the inset match
    those of the data sets in the main graph.
  }
  \label{fig:numft}
\end{figure}

A numerical computation of the conditional mean free times of the
minimal ball-piston model for a range of energy values\footnote{The
  energy values include $\ep = 1/100$, $2/100$, $\dots$, $49/100$, to
  which are added $\ep = 1/200$, $1/400$, $\dots$, $1/3200$ and $\ep =
  {1}/{2}-1/200$, ${1}/{2}-1/400$, $\dots$, ${1}/{2}-1/3200$.}
was carried out taking $\rho = (\sqrt 33 - 2)/(5 \sqrt 2) \simeq 0.5296$
and varying $\delta$ in the interval $0 < \delta \leq 0.2$. The
comparison with equation \eqref{eq:numft} is shown in
\fref{fig:numft}. All data sets collapse on the analytic result,
within an accuracy that is controlled by the number of initial
conditions. For each parameter $\delta$ and energy value $\ep$, a
number of $10^4$ initial conditions were generated with respect to the
distribution of density \eqref{eq:nuepdensity}.  A comparison between
equation \eqref{eq:taubpmft} and numerical computations of the mean
free time is shown in the inset of the same figure.

\section{Stochastic reduction and limiting Markov
  process \label{sec:stoch}} 

The mean free time \eqref{eq:taubpmft} determines the time
scale of the stochastic process of energy exchanges of the ball-piston
pair. Given an equilibrium ingoing energy configuration
$\{\eb, \ep\}$, at collision, the density of the equilibrium measure
\eqref{eq:nuepdensity} yields the probability to find the  system in
the outgoing energy configuration $\{\teb, \tep\}$.

\subsection{Stochastic kernel \label{sec:stokern}} 

Prior to resolving the collision event,  let us assume the system is
in the configuration  $(\mathbf{q},\mathbf{v})$, with position vector
$\mathbf{q}\in\dG{bp}$ and velocity vector $\mathbf{v}$ such that the
ball-piston pair has the energy configuration $\{\eb, \ep\}$,
parametrised by equation \eqref{eq:v1v2v3ep}, and such that
$\mathbf{v}\cdot\mathbf{n} < 0$. After the collision, we must
have the outgoing velocity vector $\mathbf{\tilde{v}}$,
$\mathbf{\tilde{v}}\cdot\mathbf{n} > 0$, with components

\begin{equation}
  \label{eq:tv1tv2tv3ep}
  \begin{split}
    \widetilde{v_1} & = \sqrt{2 \teb} \cos \tilde\alpha = \sigma \sqrt{2\ep}\,,\\
    \widetilde{v_2} & = \sqrt{2 \teb} \sin \tilde\alpha =\sqrt{2 \eb} \sin \alpha\,,\\
    \widetilde{v_3} &= \tilde \sigma \sqrt{2\tep} = \sqrt{2 \eb} \cos \alpha\,.
  \end{split}
\end{equation}

\noindent In particular, $\teb = \ep + \eb \sin^2 \alpha$
and $\tep = \eb \cos^2 \alpha$.

The probability per unit time of this transition is

\begin{equation}
  \frac{|\dG{bp}|}{4\pi |\Gv|}
  (\sqrt{\eb} \cos \alpha - \sigma \sqrt{\ep})_+
  \mathrm{d}\alpha\,,
  \label{eq:probalpha}
\end{equation}

\noindent which we now wish to rewrite in terms of the outgoing piston
energy $\tep$. From the third line in equation \eqref{eq:tv1tv2tv3ep},
we see that $\alpha$ can be written explicitly in terms of the ingoing and
outgoing energies,

\begin{equation}
  \sqrt{\eb} \cos \alpha = \tilde\sigma \sqrt{\tep}\,.
\end{equation}

\noindent The measure element thus transforms to

\begin{equation}
  \mathrm{d}\alpha = \mathrm{d}\tep \frac{1}{2\sqrt{\tep(\eb - \tep)}}
  \thetah(\eb - \tep) \,,
\end{equation}

\noindent where we have inserted the Heaviside step function,
$\thetah(x) = 1$ if $x \geq 0$, $0$ otherwise, to keep track of the condition
$\mathbf{v}\cdot\mathbf{n} < 0$.

Now summing over $\sigma$ and $\tilde\sigma$ and multiplying the above
expression by 2, which reflects the fact that the collision process is
independent of the sign of $v_2$, the probability per unit time
\eqref{eq:probalpha} transposes to

\begin{equation}
  \label{eq:Wdtep}
  W( \eb, \ep | \teb, \tep)\, \mathrm{d}\tep\,,
\end{equation}

\noindent where the probability density,

\begin{align}
  W( \eb, \ep | \teb, \tep)
  &
    = \frac{|\dG{bp}|}{4\pi |\Gv|}
    \sum_{\sigma, \tilde\sigma}\frac{(\tilde\sigma \sqrt{\tep} - \sigma \sqrt{\ep})_+}
    {\sqrt{\tep(\eb - \tep)}}  \thetah(\eb - \tep)
    \,,
    \cr
  &
    =
    \frac{|\dG{bp}|}{2\pi |\Gv|}
    \frac{\mathrm{max}(\sqrt{\tep}, \sqrt{\ep})}
    {\sqrt{\tep(\eb - \tep)}} \thetah(\eb - \tep)
    \,,
    \label{eq:kernel}
\end{align}

\noindent can be interpreted as the rate of probability of a transfer
of energy $\zeta = \eb - \teb$ from the ball at energy $\eb$ to the
piston whose energy changes from $\ep$ to $\tep = \ep + \zeta$, $-\ep
\leq \zeta \leq \eb$.

By construction, we recover the conditional collision rate
\eqref{eq:fmftbpeprate} after integrating equation \eqref{eq:kernel}
over $\zeta$,

\begin{equation}
  \label{eq:intkernel}
  \fmft{bp}{\ep}
  = \int \mathrm{d} \zeta\, W( \eb, \ep | \eb -  \zeta, \ep +  \zeta)
  \,.
\end{equation}

\subsection{Convergence to a Markov process \label{sec:stoconv}}  

Although equation \eqref{eq:kernel} is a property of the equilibrium
system, we argue that it also provides an accurate description of the
energy exchange process between the ball and piston away from
equilibrium, provided we consider the limiting regime of rare
interactions---that is, when the penetration length of the piston
into the domain of the ball is arbitrarily small, $\delta\ll1$.
Indeed, under this assumption, the ball and piston typically undergo
many wall-collision events between every binary collision, so that a
relaxation to equilibrium of the ball-piston pair at fixed energies
effectively takes place before the next occurrence of a binary collision.

As a result, in the limiting regime of rare interactions, the process
of energy exchanges is expected to converge to a Markov jump process
with kernel \eqref{eq:kernel}. Indeed, the Markov property essentially
means that if one knows not only the present energy partition, but
also has information about its history, this additional information
does not improve one's ability to predict the future evolution of the
energy partition. Correspondingly, convergence to a Markov process
means that if $\delta>0$ decreases, then information about the past
influences the future less and less. So we need to see that if we
start the system with a given energy partition between the disk and
piston, but away from the equilibrium measure (conditioned on the
energy surface), then, as $\delta\to0$, we measure nearly the same
jump rates as in \eqref{eq:kernel}. This is expected to hold exactly
because of the relaxation to (conditional) equilibrium during the many
wall-collisions that typically precede the first energy exchange. 

According to this argument, the time-evolution of the ball-piston pair
energy distribution $P(\{\eb, \ep\}, t)$ may be described by the
following master equation: 

\begin{align}
  \label{eq:meq}
  \partial_t P(\{\eb, \ep\}, t)
  &=
    \int \mathrm{d} \zeta\,
    \Big[
    W( \eb + \zeta, \ep - \zeta | \eb, \ep)
    P(\{ \eb + \zeta, \ep - \zeta \}, t)
    \cr
  &
    \qquad
    -     W( \eb, \ep | \eb - \zeta, \ep + \zeta)
    P(\{ \eb, \ep \}, t) \Big]\,.
\end{align}

\noindent The method used to derive this result relies on geometric
and measure-theoretic arguments. An alternative approach based on
kinetic theory, such as used in references \cite{Gaspard:2008NJP3004,
  Gaspard:200811P021, Gaspard:2009P08020}, yields the same results.

\subsection{Numerical test of the Markov property \label{sec:stonum}}  

We want to further substantiate the assertion that---with an
appropriate rescaling of time---the actual process of energy exchanges
produced by the ball-piston billiard has a well-defined limit when
$\delta\to 0$, and that the limiting process is indeed Markovian and
described by equation \eqref{eq:meq}. Based on the heuristic arguments
presented in \sref{sec:stoconv}, we should therefore check that, as
$\delta$ decreases, the collision rate and post-collision energy
distribution become arbitrarily close to the limits associated with
the process generated by equation \eqref{eq:meq}, independently of our
information about the past, that is, independently of the
(non-equilibrium) initial distribution we sample our energy partition
with. 

Moreover the notion that we have limited information about the past
reflects a lack of precise knowledge of the initial conditions. Our
initial distribution should therefore be smooth. Finally, since there
is, of course, no hope of checking that this is convergence holds for 
\emph{every} initial distribution, we pick a specific family
of smooth initial measures on the energy surface, and check
convergence for it. 

In particular we consider, for different values of $\delta$, an initial
state $(\mathbf{q},\mathbf{v})$ of the billiard map, with position
$\mathbf{q}$ uniformly distributed on the collision surface $\dG{bp}$ and
velocity $\mathbf{v}$ as in \eqref{eq:v1v2v3ep}, such that
$\mathbf{v}\cdot\mathbf{n}>0$,  with angle $\alpha$ and sign $\sigma$
now distributed \emph{away} from the distribution of density
\eqref{eq:nuepdensity}, and measure the distribution of ingoing
velocities at the first ball-piston collision event. The velocity
$\mathbf{v}$, now such that $\mathbf{v}\cdot\mathbf{n} < 0$, may again
be parametrised as in \eqref{eq:v1v2v3ep},  with values of $\alpha$
and $\sigma$ different from the outgoing initial velocity, but $\ep$
unchanged. Since we are considering a marginal velocity distribution,
apart from the two values of the sign $\sigma$, this distribution is a
function of a single real variable, $\alpha$. Irrespective of the
choice of piston energy $\ep$ (with ball energy $\eb = \tfrac{1}{2} -
\ep$), we expect to find a distribution of ingoing $\alpha$ and
$\sigma$ that, as $\delta\to0$, becomes arbitrarily close to
the distribution on the constant $\ep$ circles induced by the
equilibrium measure.

That is because, in the absence of ball-piston collisions, the
wall-collision events will typically induce a relaxation of the
billiard dynamics to the measure of density \eqref{eq:nuepdensity},
which is a true invariant measure of the non-interacting ball-piston
dynamics when their energies are fixed to the corresponding values. In
other words, when $\delta$ is small, the billiard dynamics is likely
to perform many wall collision events before first hitting the
ball-piston collision surface $\Mm{bp}$. The first hitting
distribution is thus expected to converge to the distribution of
density \eqref{eq:nuepdensity}, which happens to be an equilibrium
distribution of the ball-piston dynamics when interactions are turned
off ($\delta \equiv 0$).

To be specific, let

\begin{equation}
  \label{eq:falpha}
  h_{\ep}^{(n)} ( \alpha, \sigma) \propto (\sqrt{\eb} \cos \alpha -
  \sigma \sqrt{\ep})_+^n \,,
\end{equation}

\noindent and normalise these densities so that $\sum_{\sigma=\pm1}
\int \mathrm{d}\alpha \,   h_{\ep}^{(n)} ( \alpha, \sigma) = 1$. In
particular, the density $h_{\ep}^{(0)}$ is uniform on the set
$\mathbf{v}\cdot\mathbf{n} > 0$ and $h_{\ep}^{(1)}$ is the density
\eqref{eq:nuepdensity} induced by the equilibrium distribution, albeit
with a different normalisation.

\begin{figure}[!p]
  \centering
  \begin{subfigure}[t]{0.32\textwidth}
    \includegraphics[width=\textwidth]
    {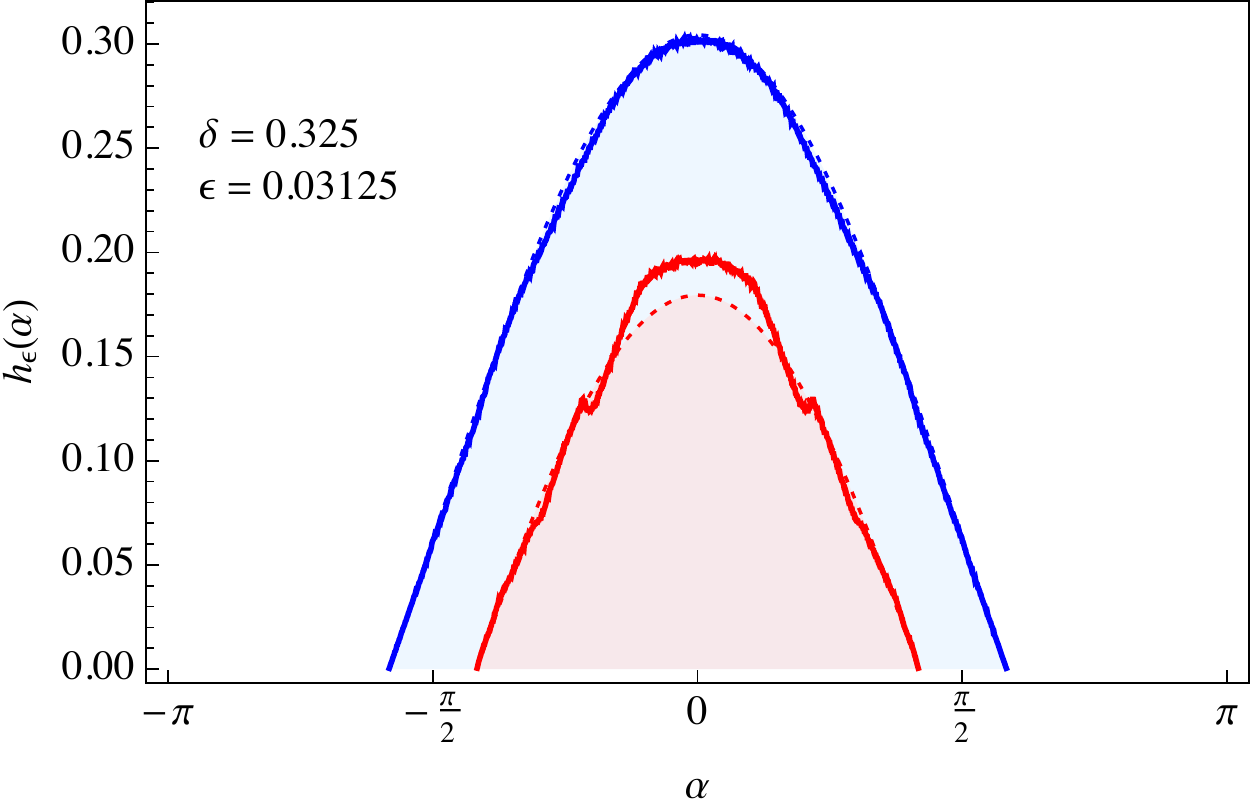}
    \caption{$\epsilon = 0.03125$, $\delta = 0.325$}
    \label{fig:rhodeltato011}
  \end{subfigure}
  \begin{subfigure}[t]{0.32\textwidth}
    \includegraphics[width=\textwidth]
    {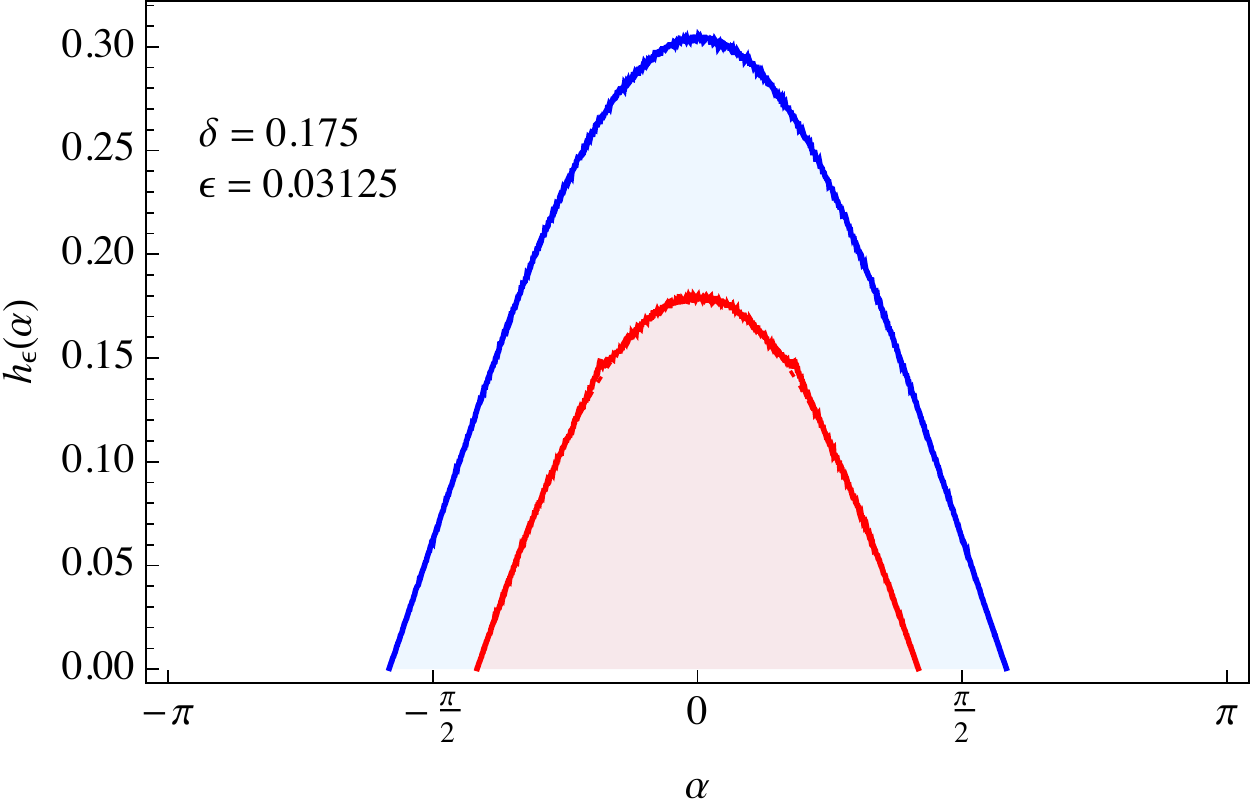}
    \caption{$\epsilon = 0.03125$, $\delta = 0.175$}
    \label{fig:rhodeltato012}
  \end{subfigure}
  \begin{subfigure}[t]{0.32\textwidth}
    \includegraphics[width=\textwidth]
    {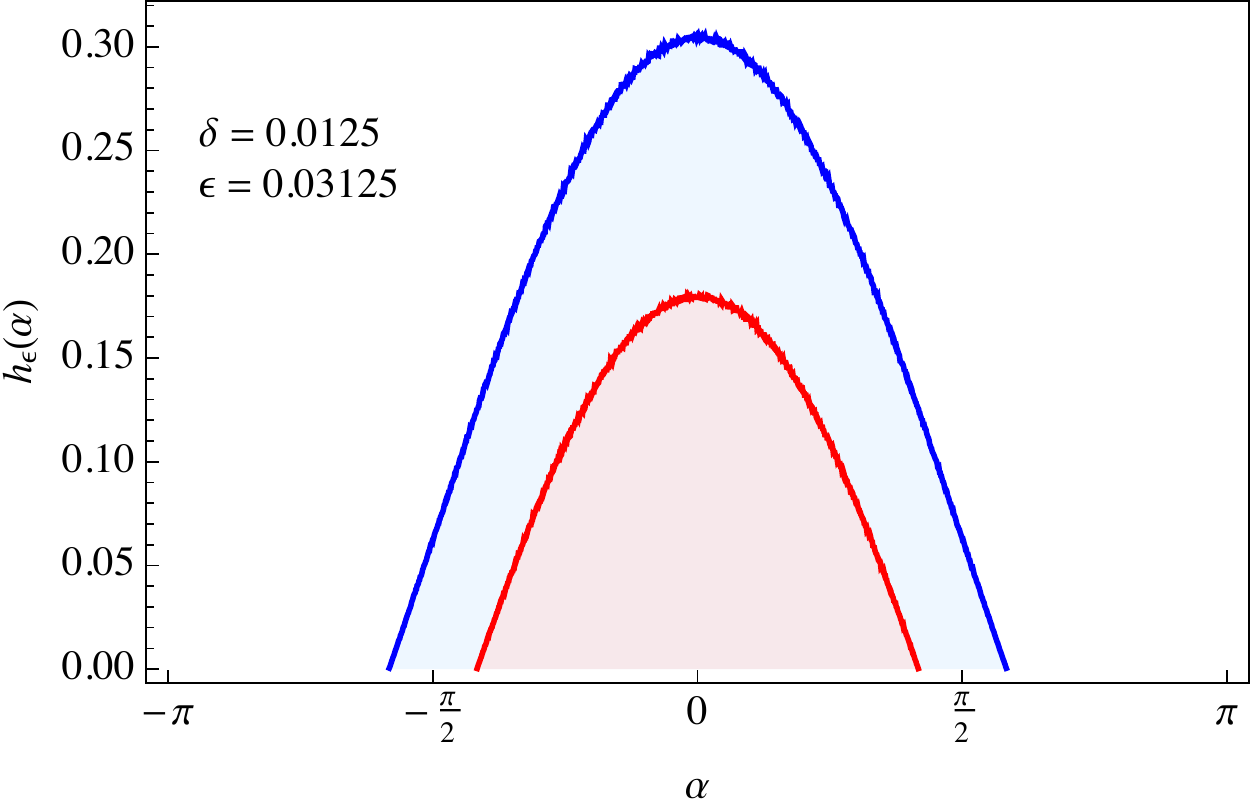}
    \caption{$\epsilon = 0.03125$, $\delta = 0.0125$}
    \label{fig:rhodeltato013}
  \end{subfigure}

  \begin{subfigure}[t]{0.32\textwidth}
    \includegraphics[width=\textwidth]
    {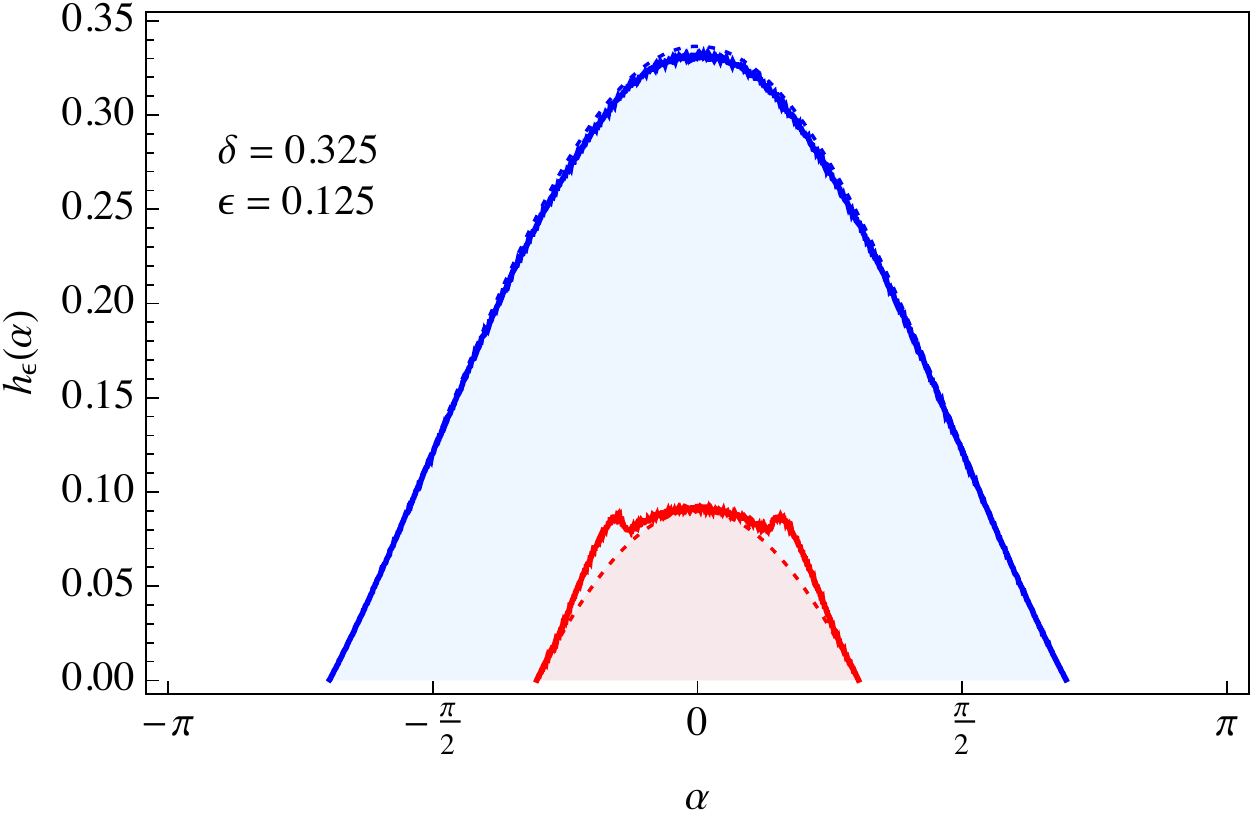}
    \caption{$\epsilon = 0.125$, $\delta = 0.325$}
    \label{fig:rhodeltato021}
  \end{subfigure}
  \begin{subfigure}[t]{0.32\textwidth}
    \includegraphics[width=\textwidth]
    {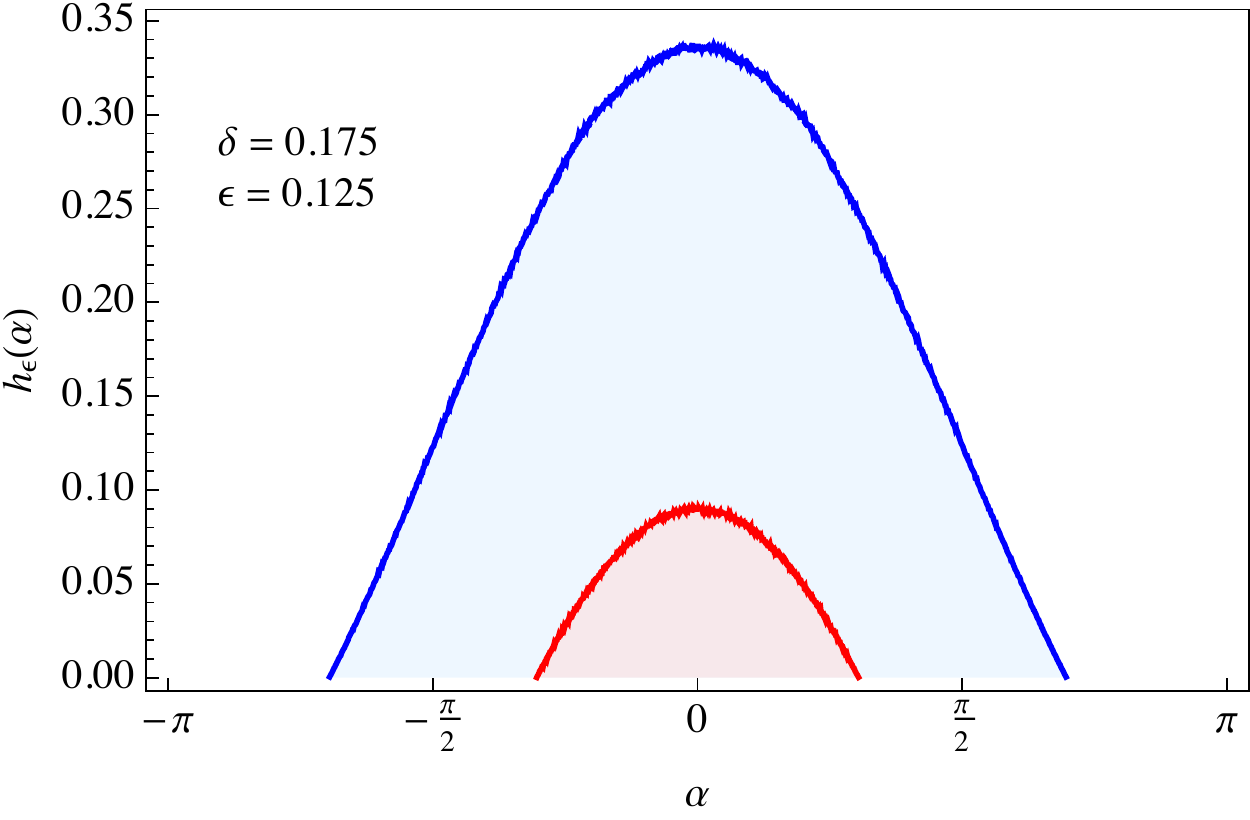}
    \caption{$\epsilon = 0.125$, $\delta = 0.175$}
    \label{fig:rhodeltato022}
  \end{subfigure}
  \begin{subfigure}[t]{0.32\textwidth}
    \includegraphics[width=\textwidth]
    {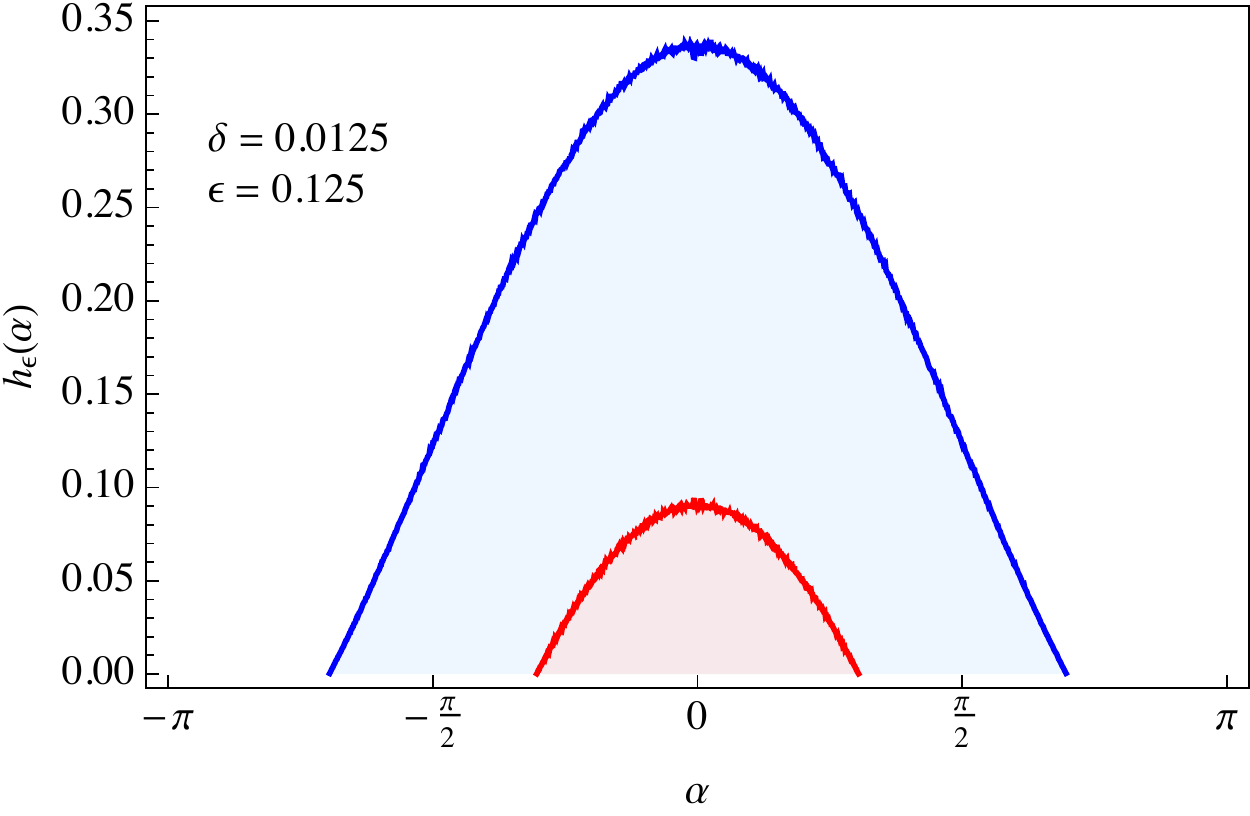}
    \caption{$\epsilon = 0.125$, $\delta = 0.0125$}
    \label{fig:rhodeltato023}
  \end{subfigure}

  \begin{subfigure}[t]{0.32\textwidth}
    \includegraphics[width=\textwidth]
    {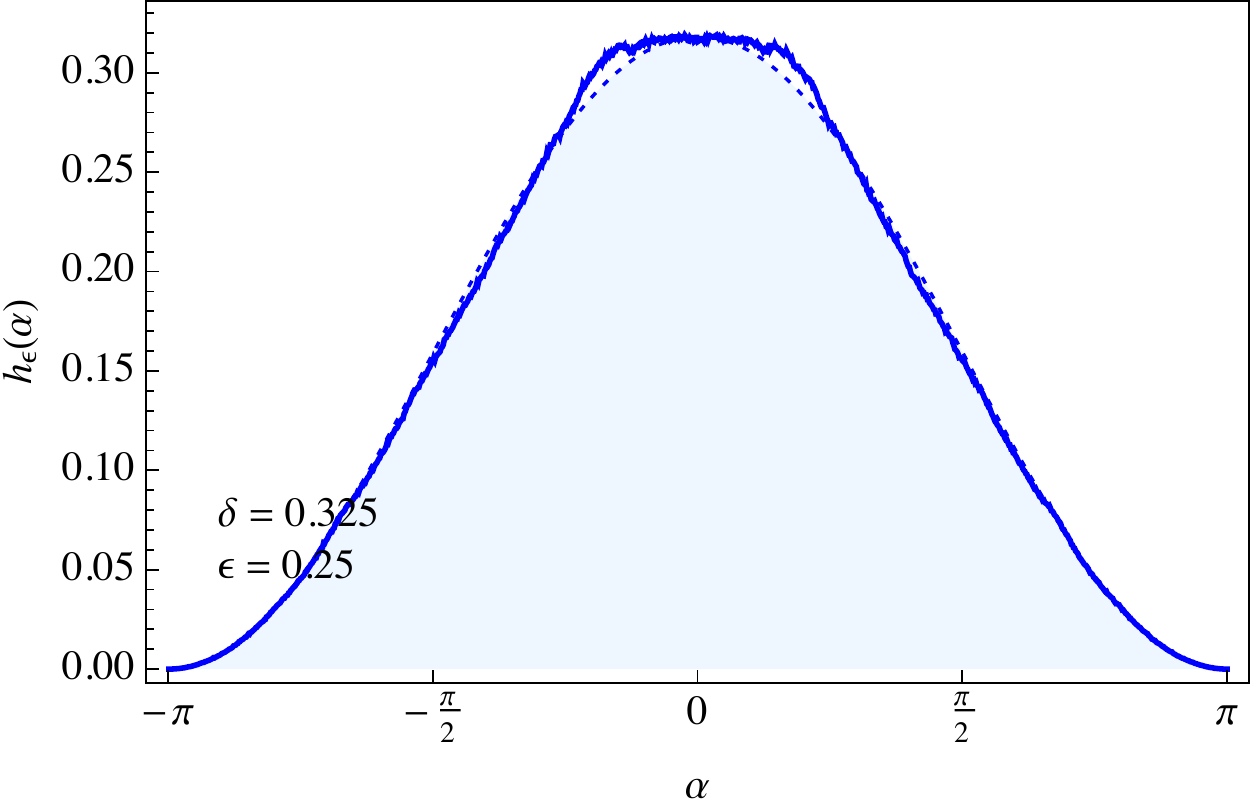}
    \caption{$\epsilon = 0.25$, $\delta = 0.325$}
    \label{fig:rhodeltato031}
  \end{subfigure}
  \begin{subfigure}[t]{0.32\textwidth}
    \includegraphics[width=\textwidth]
    {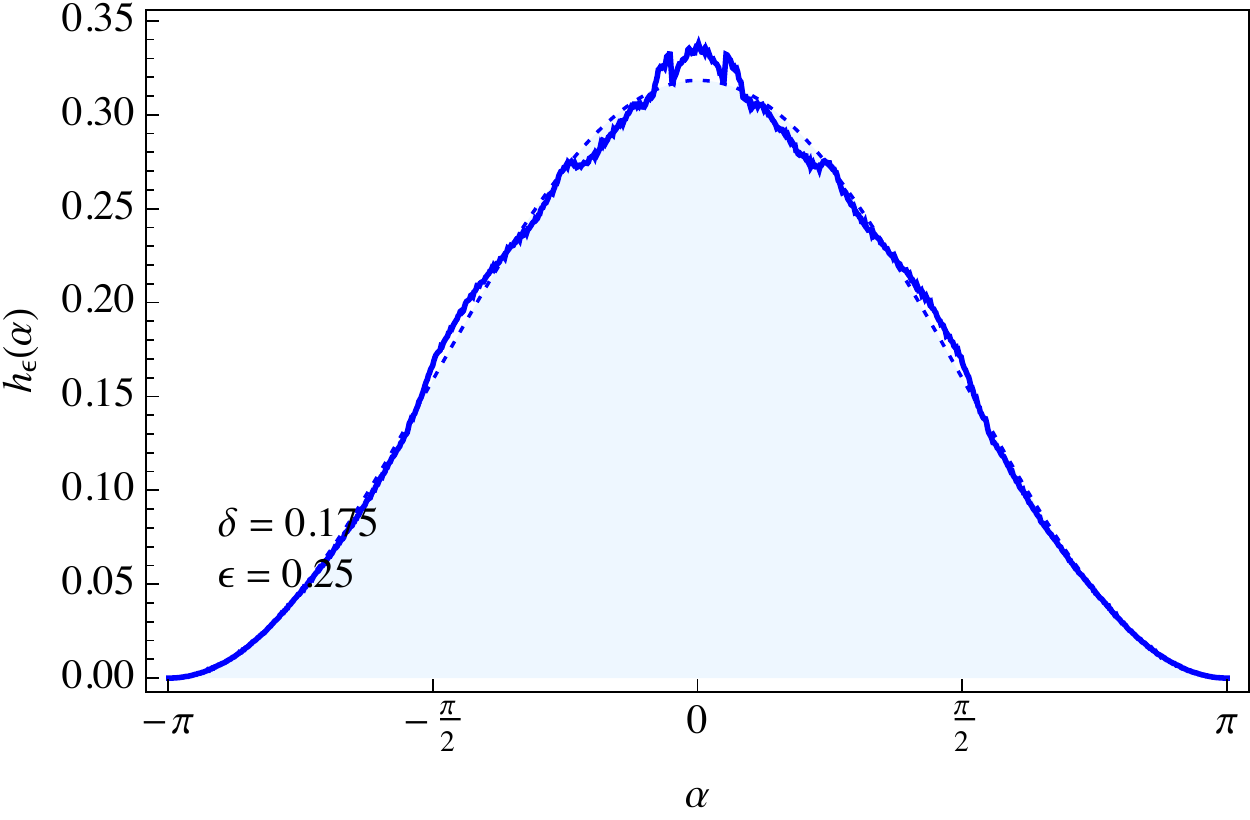}
    \caption{$\epsilon = 0.25$, $\delta = 0.175$}
    \label{fig:rhodeltato032}
  \end{subfigure}
  \begin{subfigure}[t]{0.32\textwidth}
    \includegraphics[width=\textwidth]
    {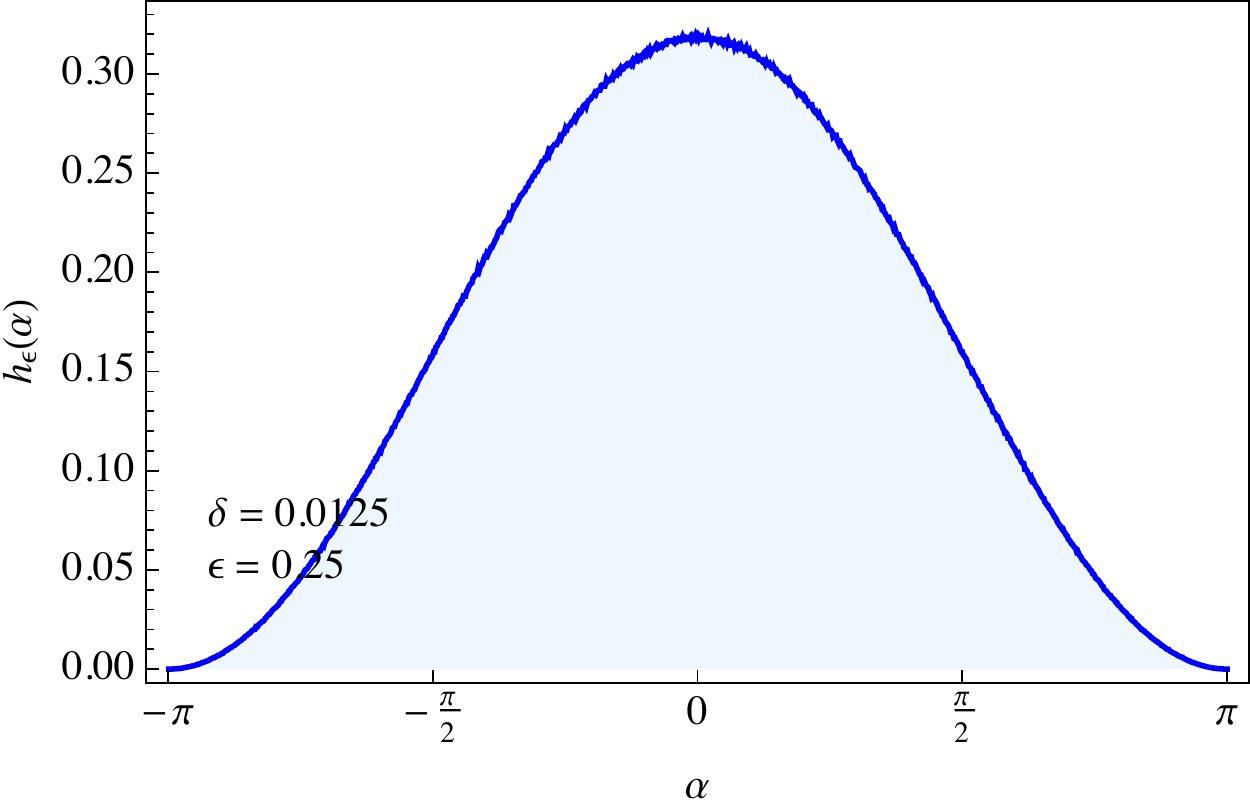}
    \caption{$\epsilon = 0.25$, $\delta = 0.0125$}
    \label{fig:rhodeltato033}
  \end{subfigure}

  \begin{subfigure}[t]{0.32\textwidth}
    \includegraphics[width=\textwidth]
    {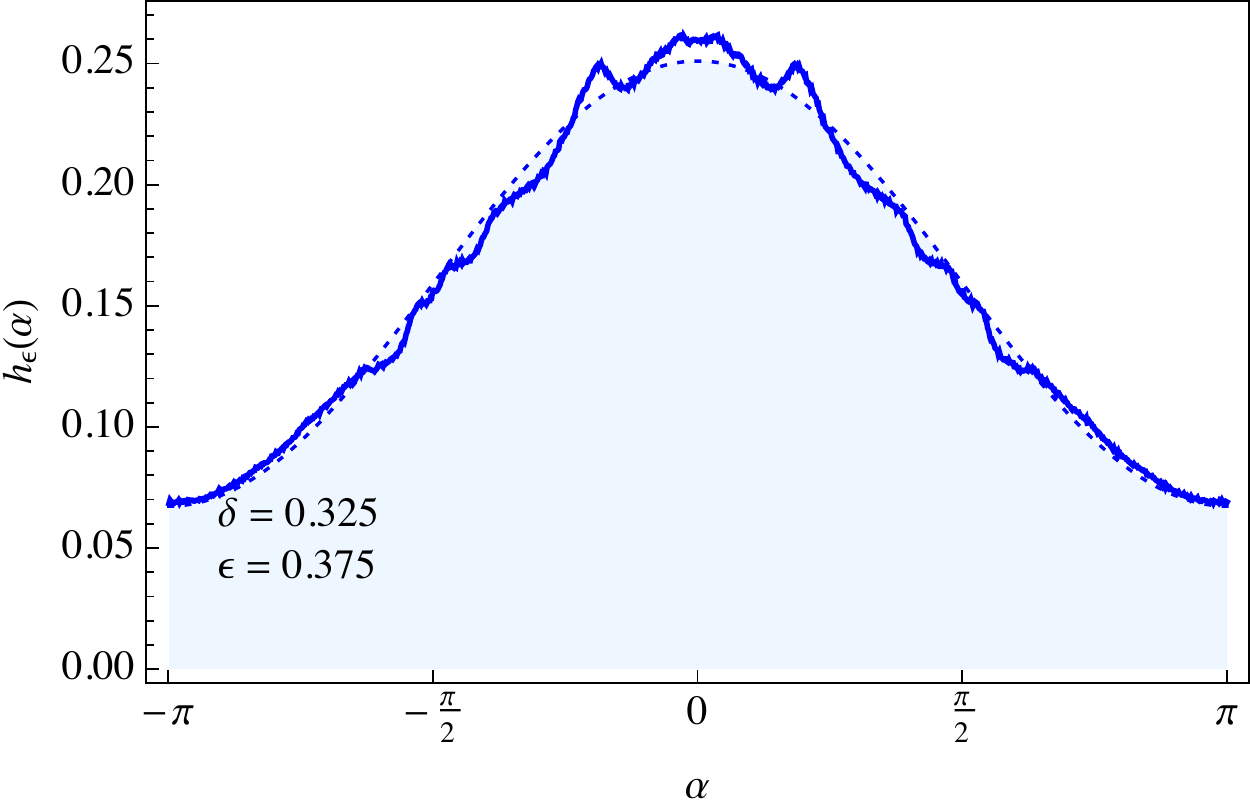}
    \caption{$\epsilon = 0.375$, $\delta = 0.325$}
    \label{fig:rhodeltato041}
  \end{subfigure}
  \begin{subfigure}[t]{0.32\textwidth}
    \includegraphics[width=\textwidth]
    {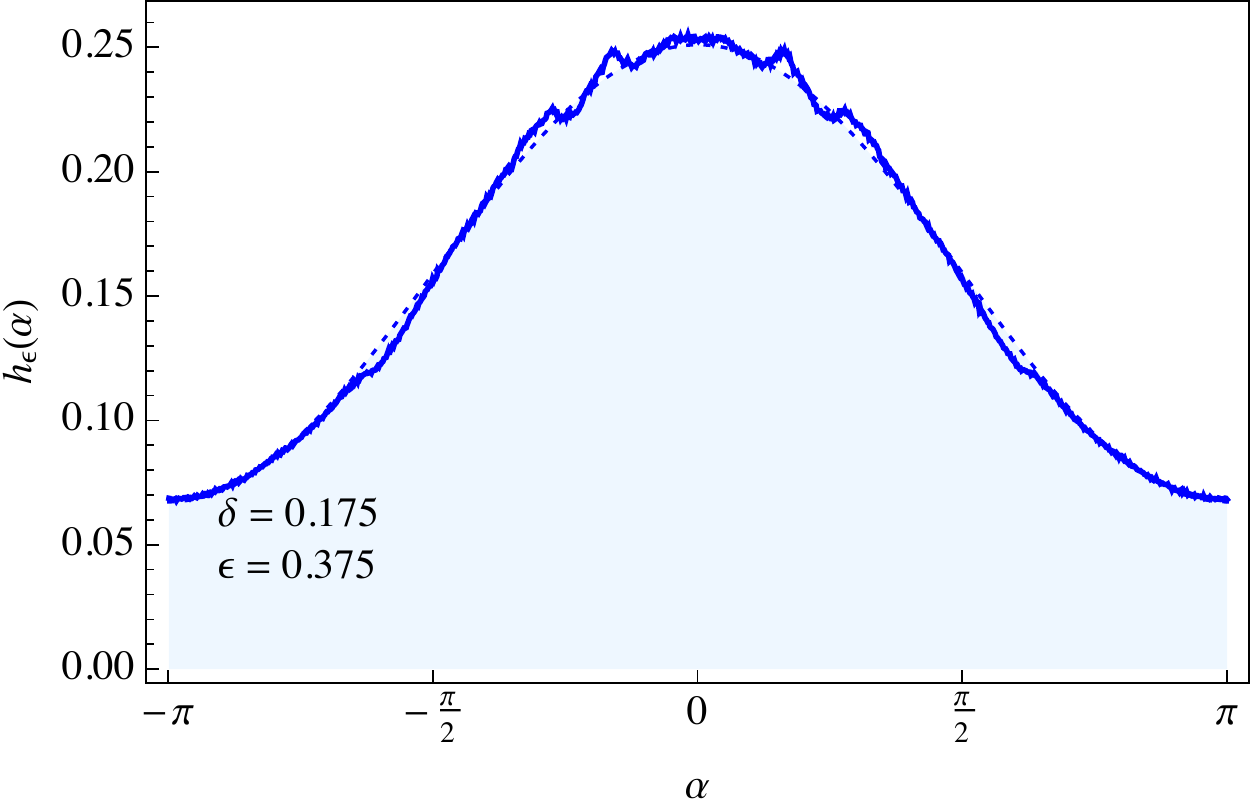}
    \caption{$\epsilon = 0.375$, $\delta = 0.175$}
    \label{fig:rhodeltato042}
  \end{subfigure}
  \begin{subfigure}[t]{0.32\textwidth}
    \includegraphics[width=\textwidth]
    {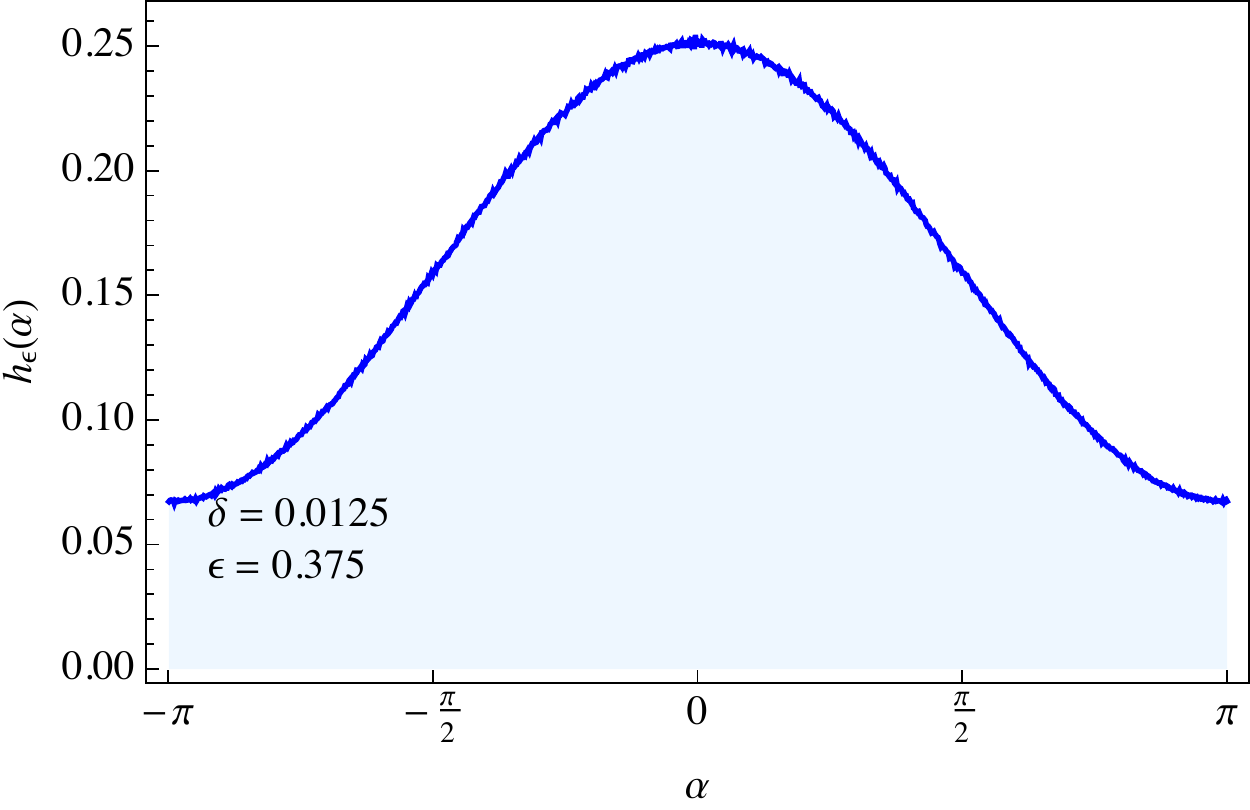}
    \caption{$\epsilon = 0.375$, $\delta = 0.0125$}
    \label{fig:rhodeltato043}
  \end{subfigure}

  \begin{subfigure}[t]{0.32\textwidth}
    \includegraphics[width=\textwidth]
    {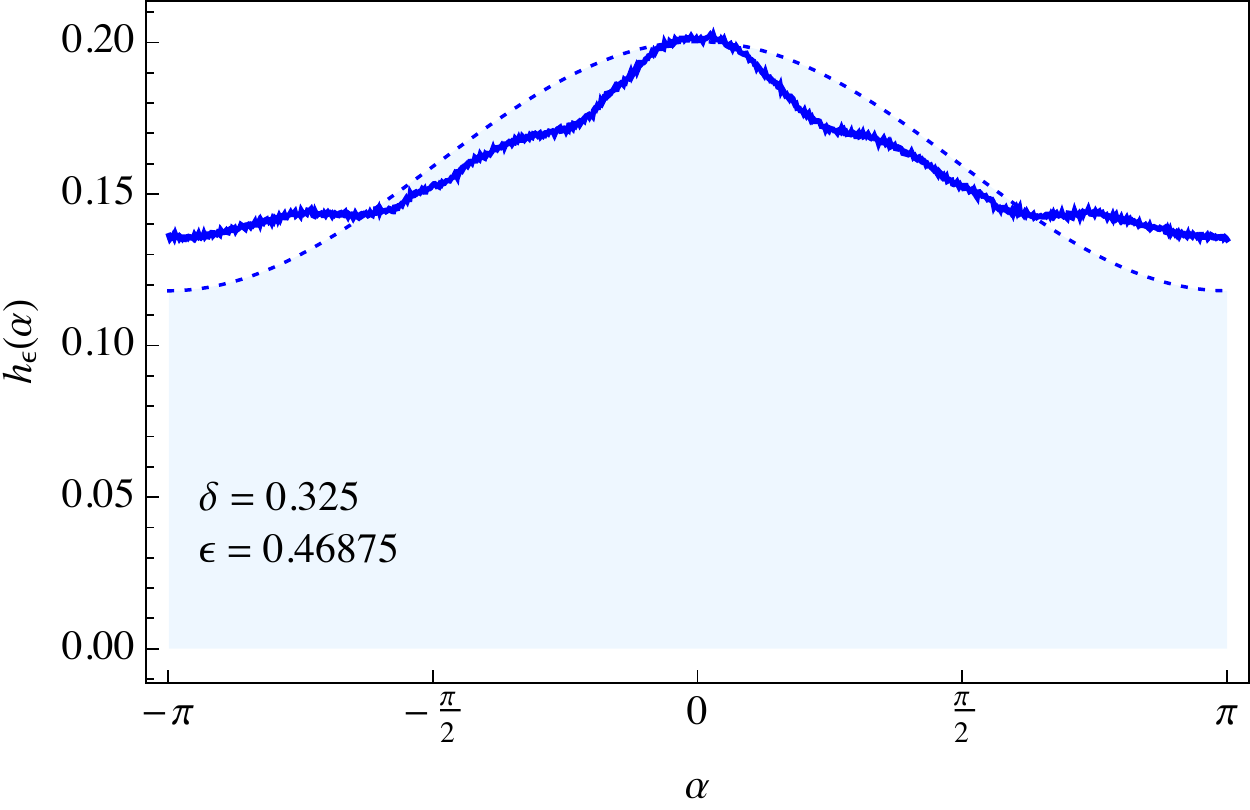}
    \caption{$\epsilon = 0.469$, $\delta = 0.325$}
    \label{fig:rhodeltato051}
  \end{subfigure}
  \begin{subfigure}[t]{0.32\textwidth}
    \includegraphics[width=\textwidth]
    {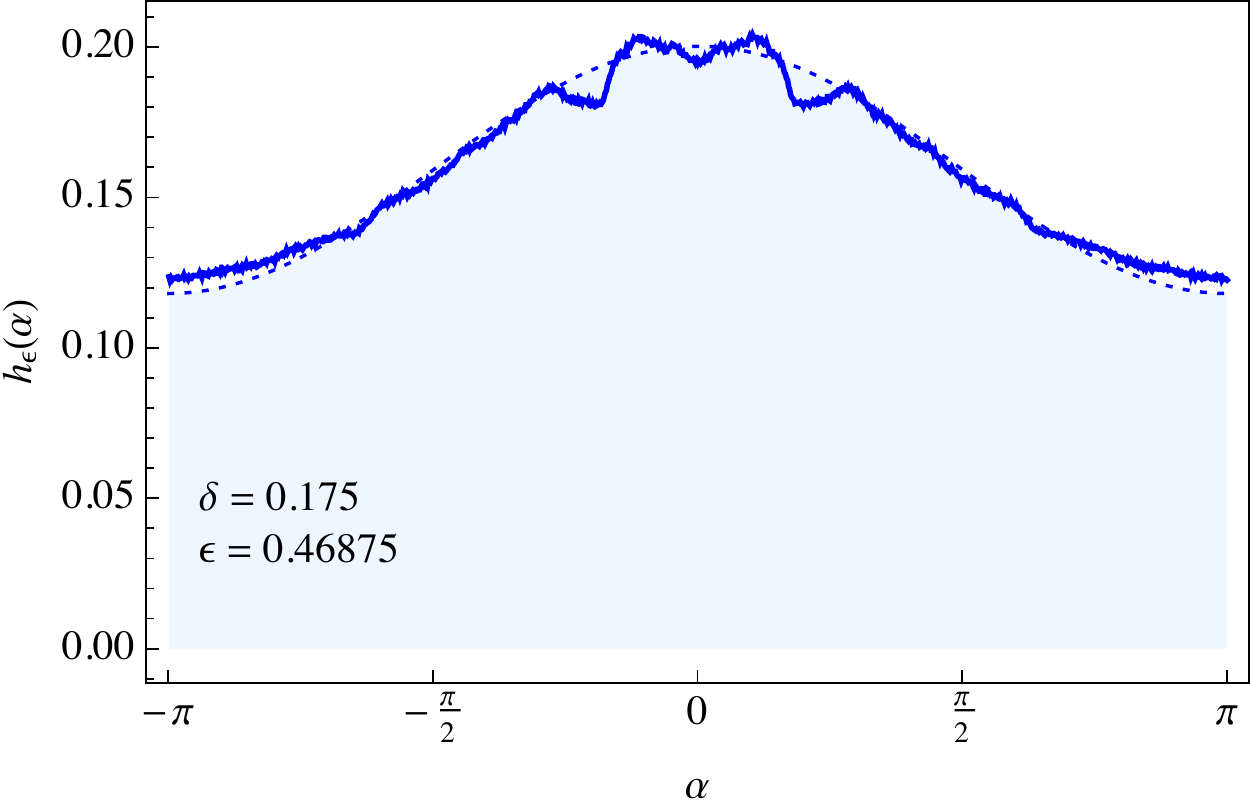}
    \caption{$\epsilon = 0.469$, $\delta = 0.175$}
    \label{fig:rhodeltato052}
  \end{subfigure}
  \begin{subfigure}[t]{0.32\textwidth}
    \includegraphics[width=\textwidth]
    {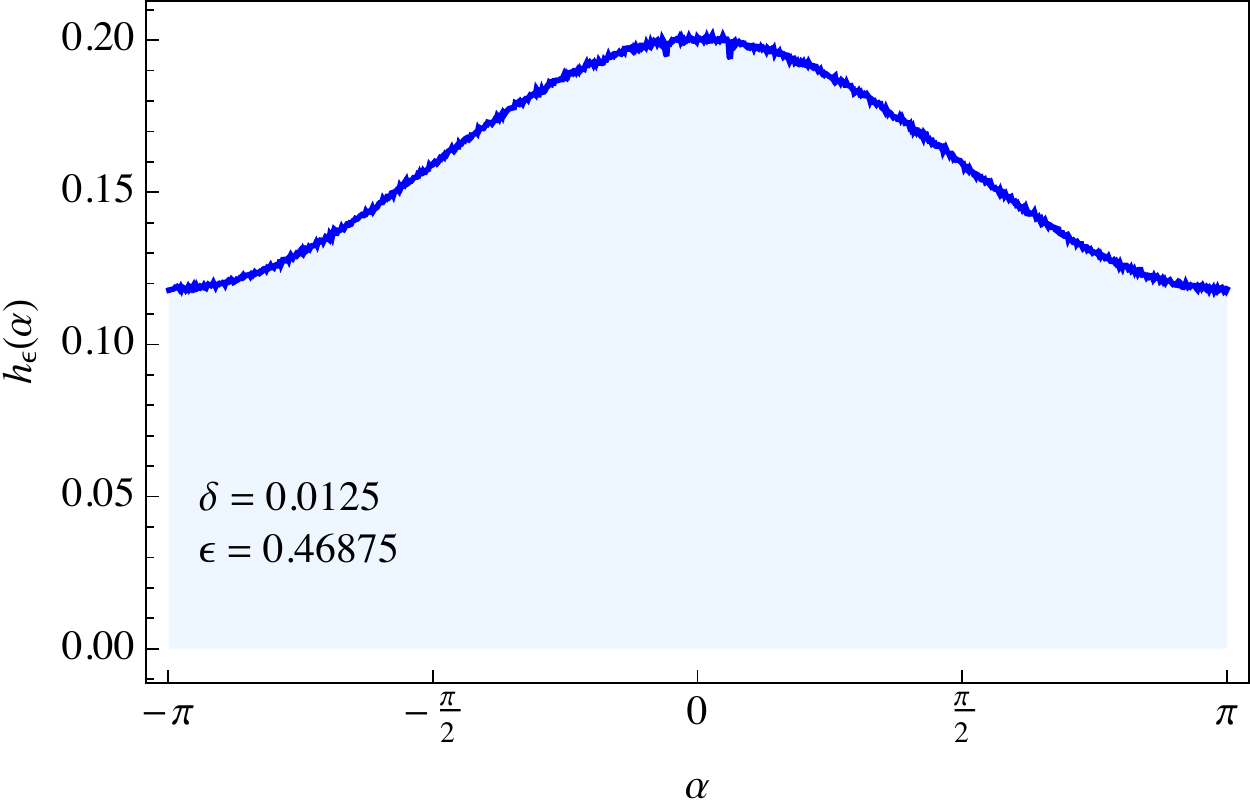}
    \caption{$\epsilon = 0.469$, $\delta = 0.0125$}
    \label{fig:rhodeltato053}
  \end{subfigure}

  \caption{Histograms of the measured ingoing first hit velocity
    distributions obtained from outgoing distributions
    \eqref{eq:falpha} with exponent $n=0$. The horizontal axes show
    the angle values $\alpha$. The parts of the densities
    corresponding to $\sigma = +1$ are shown in solid blue lines and
    to $\sigma = -1$ in solid red lines (only for $\ep < 1/4$). Every
    curve is compared to the corresponding density induced by the
    equilibrium measure (filled areas). Each row in the figure
    corresponds to a fixed value of $\ep$ and each column to a fixed
    value of $\delta$, decreasing from left to right.}
  \label{fig:falpha}
\end{figure}

\begin{figure}[htb]
  \centering
  \begin{subfigure}[t]{0.49\textwidth}
    \includegraphics[width=\textwidth]
    {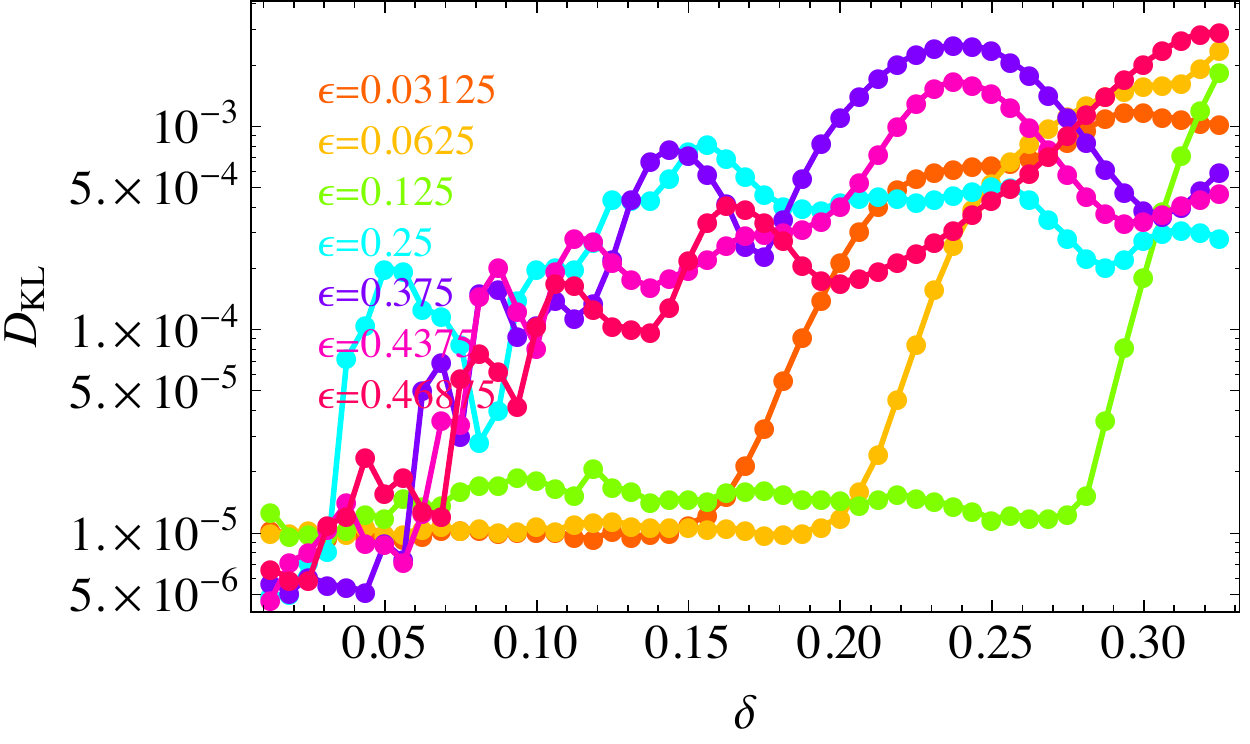}
    \caption{$n = 0$}
    \label{fig:kldn0}
  \end{subfigure}
  \begin{subfigure}[t]{0.49\textwidth}
    \includegraphics[width=\textwidth]
    {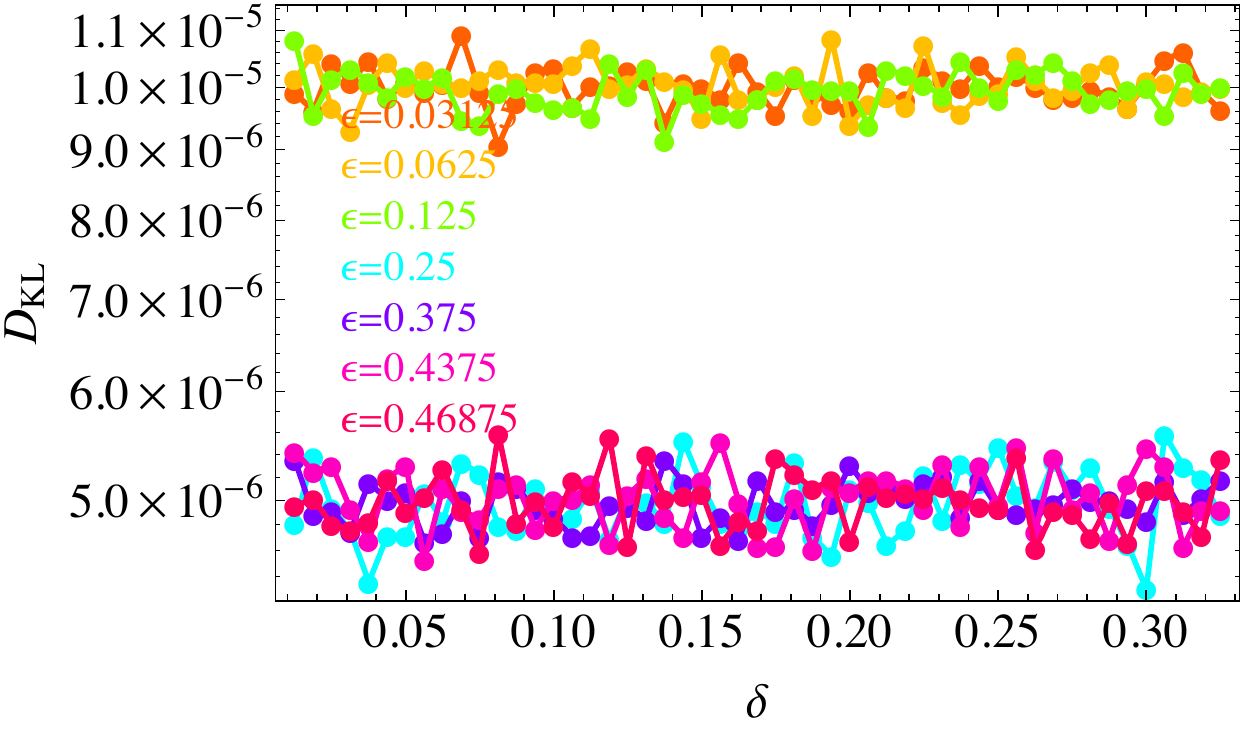}
    \caption{$n = 1$}
    \label{fig:kldn1}
  \end{subfigure}

  \begin{subfigure}[t]{0.49\textwidth}
    \includegraphics[width=\textwidth]
    {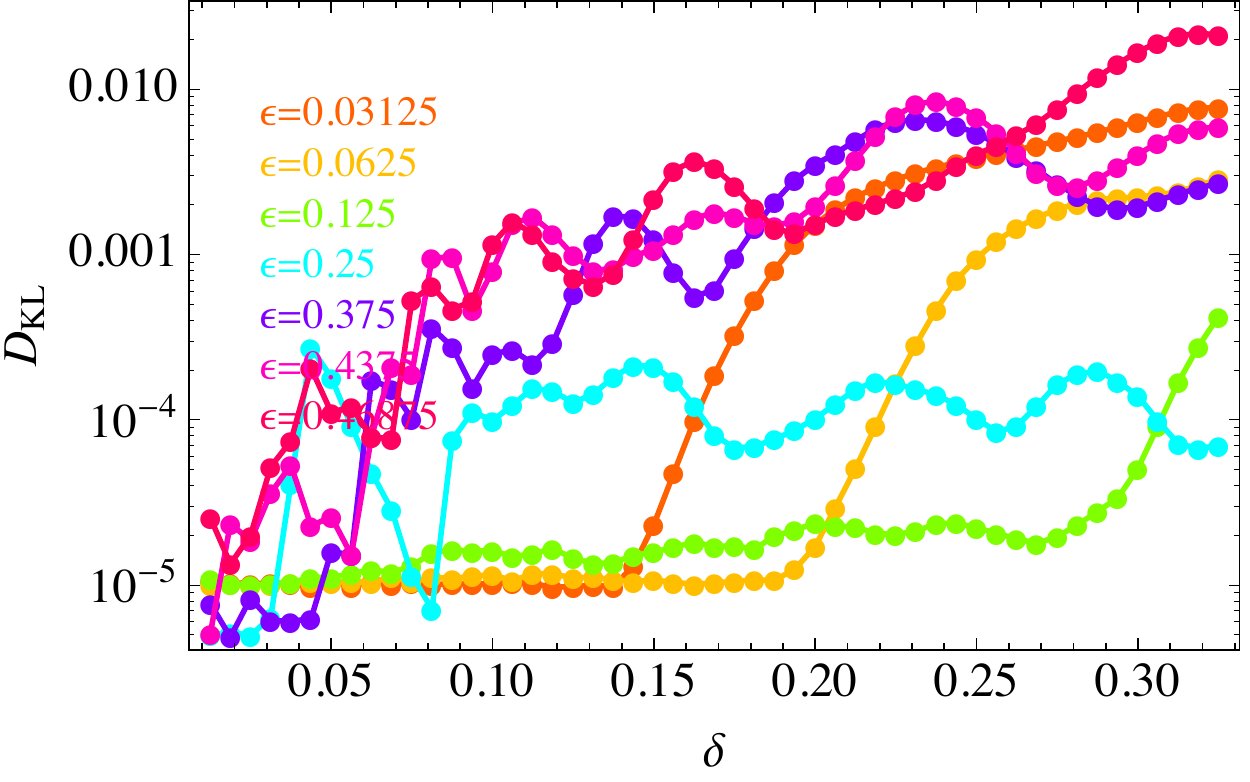}
    \caption{$n = 5$}
    \label{fig:kldn5}
  \end{subfigure}
  \begin{subfigure}[t]{0.49\textwidth}
    \includegraphics[width=\textwidth]
    {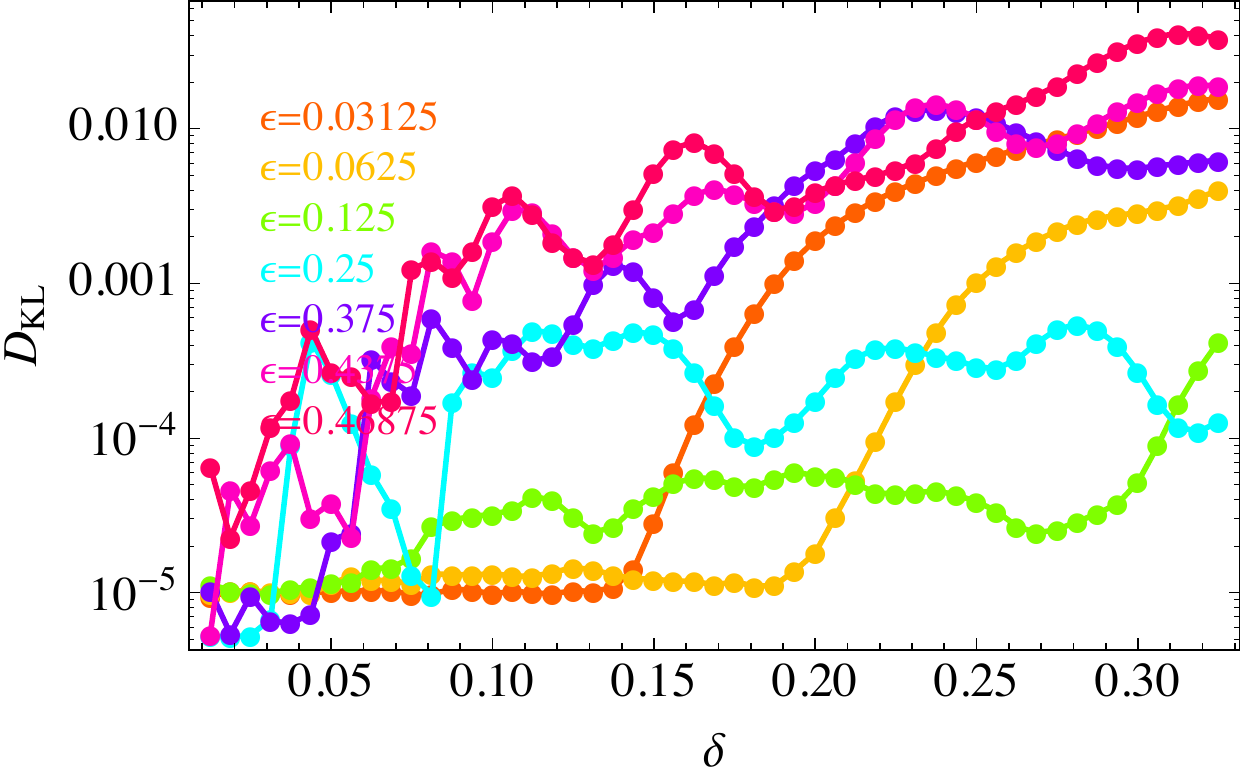}
    \caption{$n = 10$}
    \label{fig:kldn10}
  \end{subfigure}

  \caption{The Kullback-Leibler divergence \eqref{eq:kld} of
    the measured density of the ingoing velocity distributions
    relative to the equilibrium density \eqref{eq:nuepdensity}
    provides a quantitative measurement of the convergence of the
    former to the latter as the parameter $\delta\to0$. Data
    obtained by
    sampling initial conditions with respect to the densities
    $h_{\ep}^{(n)}$, (a) $n = 0$, (b) $n = 1$, (c) $n = 5$, (d) $n =
    10$ are displayed as functions of the penetration length
    $\delta$ of the piston into the ball cell. The oscillatory
    behaviour observed for some values of the energies appears to be
    logarithmic with respect to $\delta$. In panel (b), we
    compare the numerically obtained values
    of the equilibrium distribution $h_{\ep}^{(1)}$ to its analytic
    expression \eqref{eq:nuepdensity}, and thus obtain a useful
    benchmark  to gauge the accuracy within which equilibrium
    statistics can be reached. The
    factor $2$ between the sets with $\ep < 1/4$ and $\ep \ge 1/4$
    is due to the fact that the phase space is divided into twice as
    many cells when $\ep < 1/4$ compared to $\ep \ge 1/4$.
  }
  \label{fig:kld}
\end{figure}

In \fref{fig:falpha}, we plot the histograms of the ingoing velocity
distributions obtained by sampling initial conditions with respect to
the non-equilibrium density $h_{\ep}^{(0)}$. Each subfigure
corresponds to a different value of $\delta$, varied horizontally, and
$\ep$, varied vertically. Histograms are measured by dividing the
intervals of allowed values of $\alpha$ into $10^3$ bins. As seen from
the figure, the differences between the measured distributions and the
corresponding distributions induced by the equilibrium measure are
never large, but are most noticeable when $\delta$ is large. To
quantify the convergence of the measured distributions to those
induced by the equilibrium measure, we computed the (coarse grained)
relative entropy of the measured distribution with respect to
$h_{\ep}^{(1)}$, also known as  Kullback-Leibler divergence
\cite{Kullback:1951information}. Using the notation $\tilde{h}_{\ep}^{(n)}$
to denote the ingoing distribution into which the initial distribution
$h_{\ep}^{(n)}$ evolves until the first energy exchange,

\begin{equation}
  \label{eq:kld}
  D_\textsc{kl}(\tilde{h}_{\ep}^{(n)}|h_{\ep}^{(1)}) = \sum_{\sigma = \pm1} \int
  \mathrm{d}\alpha\,
  \tilde{h}_{\ep}^{(n)}(\alpha, \sigma) \log
  \frac {\tilde{h}_{\ep}^{(n)}(\alpha, \sigma)}
  {h_{\ep}^{(1)}(\alpha, \sigma)}\,,
\end{equation}

\noindent where the integral over $\alpha$ is evaluated by summing
the measured averaged density over the total number of bins. The
results of measurements of this quantity using different
outgoing velocity distributions $h_{\ep}^{(n)}$, $n = 0, 1, 5, 10$, and
exact values for the density $h_{\ep}^{(1)}(\alpha, \sigma)$ in the
denominator of equation \eqref{eq:kld} are shown in
\fref{fig:kld}. Whereas the decay to the equilibrium noise level of
the Kullback-Leibler divergence with the parameter $\delta$
appears to be qualitatively different when the piston energies are
larger than the ball energies or vice versa,  our measurements clearly
show a systematic return to the statistics induced by the equilibrium
measure as the parameter $\delta\to0$ and thus provide a confirmation
of the observations drawn from \fref{fig:falpha}.

\subsection{Moments of the kernel}

We end this section by noting that the stochastic evolution
\eqref{eq:meq} proves particularly useful to study energy exchanges in
rarely interacting systems consisting of many particles, such as the
ball-piston gas shown in \fref{fig:bpmanycells}. In this context, we
note that the first three moments of the energy transfer rate share the
symmetries observed in other models \cite{Gaspard:2009P08020}.

Thus, given the canonical ball-piston energy distribution, which is
the product of Gamma distributions of shape parameters respectively
$\tfrac{1}{2}$ and $1$, and common scale parameter (the temperature)
$\beta^{-1}$,

\begin{equation}
  \label{eq:canmeasure}
  \Pceq(\eb, \ep)
  = \frac{\beta^{3/2}}{\sqrt{\pi \ep}} \exp[- \beta (\eb + \ep) ]\,,
\end{equation}

\noindent the zeroth moment of the energy transfer rate, similar to
$\fmft{bp}{\ep}$,  equation \eqref{eq:intkernel}, but without the
assumption $\eb = \tfrac{1}{2} - \ep$,

\begin{align}
  f(\eb, \ep)
  &\equiv \int \mathrm{d} \zeta\,
    W( \eb, \ep | \eb - \zeta, \ep +  \zeta)
  \,,
  \label{eq:intkernelnuf}
  \\
  &=
    \frac{|\dG{bp}|}
    {|\Gv|}
    \begin{cases}
    \frac{1} {\pi}
    \left[
      \sqrt{\eb - \ep} + \sqrt{\ep} \arcsin
      \sqrt{\frac{\displaystyle \ep}{\displaystyle \eb}}
    \right]  \,,
    &  \eb > \ep\,,\\
    \tfrac{1}{2} \sqrt{\ep}\,,
    & \eb \leq \ep \,,
  \end{cases}
      \nonumber
\end{align}

\noindent the first moment,

\begin{align}
  j(\eb, \ep)
  &\equiv \int \mathrm{d} \zeta\, \zeta \,
    W( \eb, \ep | \eb - \zeta, \ep +  \zeta)
    \,,
    \label{eq:intkerneleta}
  \\
  &=
    \frac{|\dG{bp}|}
    {|\Gv|}
    \begin{cases}
      \frac{1}{6\pi}
      \Big[
      (4 \eb - 7 \ep) \sqrt{\eb - \ep}
      \\
      \quad + 3 ( \eb - 2 \ep) \sqrt{\ep}
      \arcsin \sqrt{\ep/\eb}
      \Big]\,,
      & \eb > \ep\,,\\
      \frac{1}{4} (\eb - 2 \ep) \sqrt{\ep} \,,
      & \eb \leq \ep\,,
    \end{cases}
        \nonumber
\end{align}

\noindent and the second moment,

\begin{align}
  h(\eb, \ep)
  &\equiv
    \int \mathrm{d} \zeta\, \zeta^2
    W( \eb, \ep | \eb - \zeta, \ep +  \zeta)
    \,,
    \label{eq:intkerneleta2}
  \\
  &=
    \frac{|\dG{bp}|}
    {|\Gv|}
    \begin{cases}
      \frac{1}{15\pi}
      \Big\{
      8(\eb - \ep)^{5/2} + 15
      \Big[
      -\frac{3}{8} \ep ( \eb - 2 \ep) \sqrt{\eb - \ep}
      \\
      \quad + \sqrt{\ep}
      ( \frac{3}{8} \eb^2 - \eb \ep + \ep^2 )
      \arcsin \sqrt{\ep/\eb}
      \Big]
      \Big\}\,,
      & \eb > \ep\,,\\
      \tfrac{1}{2}\sqrt{\ep}       ( \frac{3}{8} \eb^2 - \eb \ep + \ep^2 )\,,
      & \eb \leq \ep\,,
    \end{cases}
        \nonumber
\end{align}

\noindent all satisfy the following identities, involving averages
with respect to the canonical measure \eqref{eq:canmeasure}:

\begin{equation}
  \label{eq:nucanavg}
  \canavg{f(\eb, \ep)}{\beta}
  = \frac{\beta^2}{2} \canavg{(\eb - \ep) j(\eb, \ep)}{\beta}
  = \frac{\beta^2}{2} \canavg{h (\eb, \ep)}{\beta}
  = \frac{1}{\sqrt{2 \pi \beta}}
  \frac{|\dG{bp}|}
  {|\Gv|}
  \,.
\end{equation}

In the limit $\delta\to 0$ of rare interactions, this is

\begin{equation}
  \label{eq:nucanavgrare}
  \lim_{\delta\to0}
  \delta^{-2}  \canavg{f(\eb, \ep)}{\beta}
  =
  \frac{1}{\sqrt{\pi \beta}}
  \Big[
  1 - \lambda - \rho^2 \big( \pi - 4 \arctan \lambda \big)
  \Big]^{-1}\,,
\end{equation}

\noindent which provides an approximation of the heat conductivity of
the ball-piston gas. This will be the subject of a separate publication.

\section{Conclusions \label{sec:conc}}

In reference \cite{Gaspard:2008PRL101} a family of billiard models was
introduced in the hope of proving suitable for deriving the heat
equation. It is moreover believed it will be possible to determine the
actual expression of the associated coefficient of heat conductivity.
Such an achievement would bring to completion a programme aiming at
explaining macroscopic laws from deterministic microscopic
assumptions, one of mathematical physics great outstanding challenges.

These models, which combine the kinetics of gases of hard balls with
the periodic structure of crystalline solids, lend themselves to a
systematic analysis, whose tools were made available in no small part
thanks to the pioneering works of David Ruelle and Yasha Sinai. The
authors of \cite{Gaspard:2008PRL101} outlined a simple two-step
strategy to attain their goal: (i) going from the microscopic scale to
a mesoscopic one (micro-to-meso), and (ii) from that scale to the
macroscopic one (meso-to-macro).  Moreover, they also realised their
programme on the level of analytic calculations with precise physical
meaning.

Among realistic models to study Fourier's law, billiard models are
generally most amenable to a rigorous derivation of both
mesoscopic and macroscopic laws from deterministic microscopic
assumptions, however delicate their technical analysis. It has
therefore been a top priority of the community
to provide a mathematically sound proof completing the approach
outlined above. Our main goal here has been to suggest a new addition
to the family of Gaspard-Gilbert models, for which a mathematical
treatment of the micro-to-meso step is a distinctly realistic task.

In this paper, we focused on the description of the ball-piston model
and the computation of several quantities characterising its
statistical properties. The limit of rare ball-piston interactions
provides a meaningful interpretation, both physically and
mathematically, of some of these properties at the level of a
mesoscopic description. Namely, energy exchanges are described by a
Markov jump process with a precise form of the transition kernel.
The complete mathematical discussion is postponed to subsequent
publications (as to the first of these, see \cite{balint:2015winp}).

We have, in addition, devised a statistical procedure, making use of
the Kullback-Leibler divergence, to test quantitatively whether the limit
of rare interactions of the minimal ball-piston model indeed possesses
the Markov property. Our numerical results are affirmative and help
shed new light on the approach to this limit. Our procedure can be put
to use in other models of the GG family and its application will be
described elsewhere.

\appendix
\numberwithin{equation}{section}
\normalsize

\section{Collision volume \label{app:vol}}

To determine the volume $|\Gv|$ of configuration space, note that
when $q_3 > (1 - \lambda)/2$, the volume of all possible positions $q_1$
and $q_2$ is

\begin{equation}
  \label{eq:areax1x2}
  4 \int_0^{(1 - \lambda)/2} \mathrm{d}\, q \left[
    \tfrac{1}{2} - \sqrt{\rho^2 -
      \left(q - \tfrac{1}{2} \right)^2} \right]
  = 1 - \lambda - \rho^2 \big( \pi - 4 \arctan \lambda \big)
  \,.
\end{equation}

\noindent When $q_3 < (1 - \lambda)/2$, we must subtract from the
above area the quantity

\begin{align}
  2 & \int_{q_3}^{(1-\lambda)/2}
      \mathrm{d}\, q
      \left[ \tfrac{1}{2} - \sqrt{\rho^2 -
      \left(q - \tfrac{1}{2} \right)^2} \right]
      = \frac{1}{2} - \frac{\lambda}{4} - q_3
      - \frac{1 - 2 q_3}{4} \sqrt{4 \rho^2 - (1 - 2 q_3)^2}
      \cr
  &\quad + \rho^2 \left[
    \arctan \lambda -
    \arctan \frac{1 - 2 q_3}
    {\sqrt{4 \rho^2 - (1 - 2 q_3)^2}}
    \right]
    \,.
    \label{eq:area1}
\end{align}

\noindent We therefore obtain the total volume of configuration space
\eqref{eq:Gv}  by multiplying equation \eqref{eq:areax1x2} by $\lambda
+ 2 \delta$ and subtracting the integral of equation \eqref{eq:area1}
over $q_3$ from $(1 - \lambda)/2 - \delta$ to $(1 - \lambda)/2$.

To compute the collision surface integral, note that the position
coordinates on $\dG{bp}$ are bounded according to

\begin{equation}
  \label{eq:dGc1}
  \begin{split}
    & \tfrac{1}{2}(1 - \lambda)  - \delta \leq q_1 = q_3 \leq
    \tfrac{1}{2} (1 - \lambda)  \,,
    \\
    &-\tfrac{1}{2} + \sqrt{\rho^2 -
      (q_1-\tfrac{1}{2})^2} \leq q_2 \leq \tfrac{1}{2} - \sqrt{\rho^2 - (q_1 - \tfrac{1}{2})^2}
  \end{split}
\end{equation}

\noindent Its projection on the $(q_1, q_2)$ plane is the area
\eqref{eq:area1} evaluated at $q_3 = (1 - \lambda)/2  - \delta$. Since
the surface itself makes an angle $\pi/4$ with respect to the  $(q_1,
q_2)$ plane, we obtain for $|\dG{bp}|$ the expression given by equation
\eqref{eq:dGc}.

\section{Ball-wall and piston-wall return times \label{app:wall}}

Wall collision return times of the piston and ball can be computed
in ways similar to equation \eqref{eq:taubpmft},

\begin{equation}
  \label{eq:tauwall}
  \begin{split}
    \taumft{pw} &= \frac{4 |\Gv|}{|\dG{pw}|}\,,
    \\
    \taumft{bw} &= \frac{4 |\Gv|}{|\dG{bw}|}\,,
  \end{split}
\end{equation}

\noindent where $|\dG{pw}|$ and  $|\dG{bw}|$ are the areas of piston-wall and
ball-wall collisions. The former corresponds to the area of all
positions $q_1$ and $q_2$ such that $q_3 = ( 1 \pm \lambda )/2 \pm
\delta$, which is parallel to the $(q_1, q_2)$ plane and twice the area
\eqref{eq:areax1x2} minus the projection of the collision surface
$|\dG{bp}|$ \eqref{eq:dGc} on this plane, and the latter to the positions
$q_1, q_2$ such that $(q_1 \pm \tfrac{1}{2})^2 + (q_2 \pm
\tfrac{1}{2})^2 = \rho^2$ while $q_1 < q_3$, with $q_3$ integrated over
the interval \eqref{eq:x1x2x3}. That is,

\begin{align}
  |\dG{pw}|
  &=
    2\Big [1 - \lambda - \rho^2 \big( \pi - 4 \arctan \lambda \big) \Big]
    - \frac{1}{\sqrt2}|\dG{bp}|\,,
\end{align}

\noindent and

\begin{align}
  |\dG{bw}|
  & =
    \rho ( \lambda + 2 \delta )
    \Big(
    8 \arcsin \frac{1 - \lambda}{2 \sqrt 2 \rho}
    - \arcsin \frac{\lambda + 2 \delta}{2 \rho}
    + \arcsin \frac{\lambda}{2 \rho}
    \Big)
    \cr
  &
    \quad + \rho \big [1 -
    \sqrt{1 - 4 \delta ( \lambda + \delta) }
    \big]
    \,.
  \label{eq:dGb}
\end{align}

\section{Restriction of the invariant measure on $\Mm{bp}$ to fixed
  energy configurations \label{app:cnuep}} 

\begin{figure}[hbt]
  \centering
  \includegraphics[width=0.65\textwidth]
  {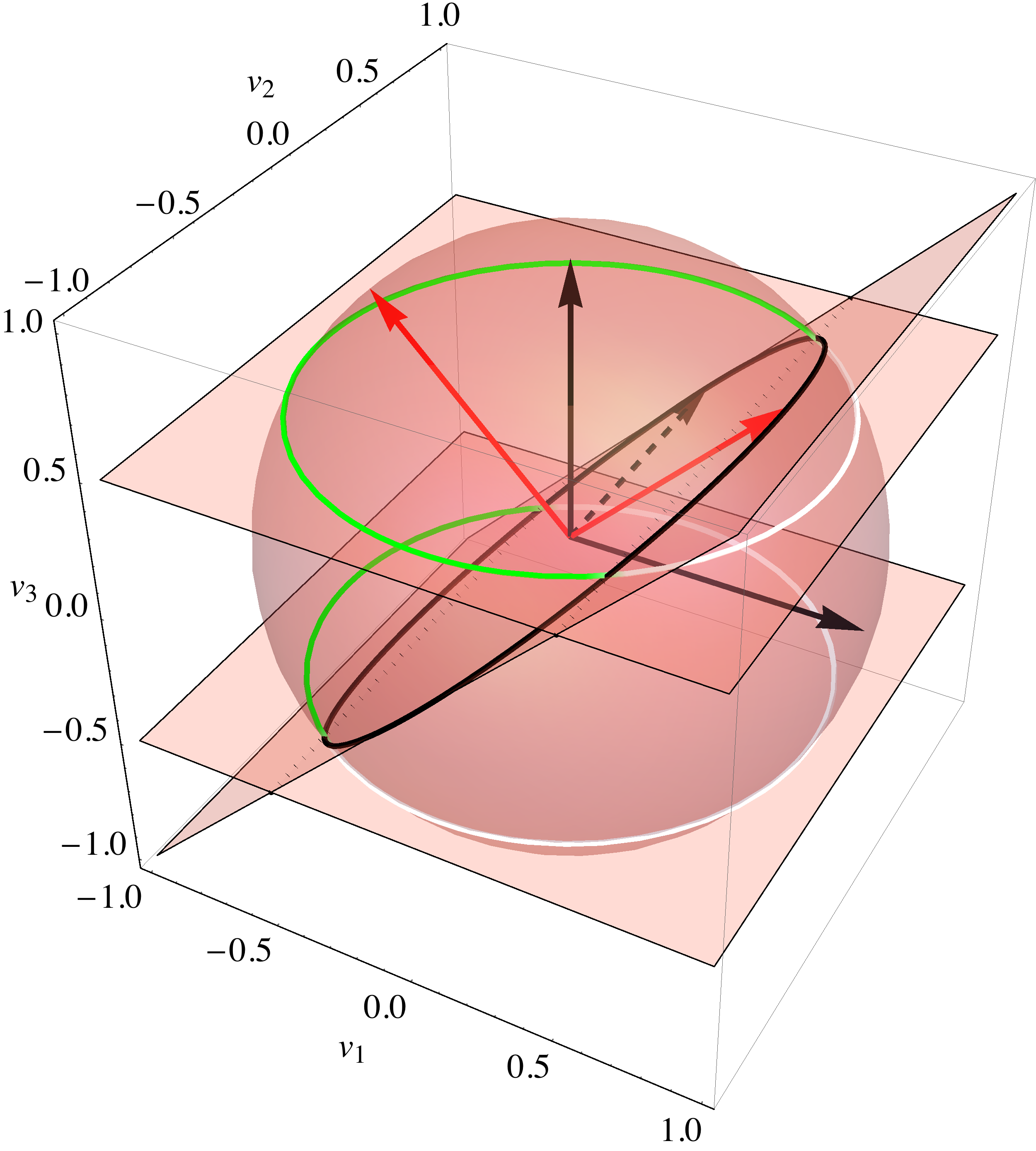}
  \caption{Equation \eqref{eq:cnuepvpart} has either one ($1/4
    \leq \ep \leq \tfrac{1}{2}$) or two ($0 \leq \ep < 1/4$) contributions,
    given by the integrals along the arc-circles on the
    hemisphere of velocity coordinates whose projection along the
    normal to the collision surface \eqref{eq:n} is positive. Here
    $\ep = 1/8$ and the two arc-circles at $v_3 = \pm \tfrac{1}{2}$
    contributing to equation \eqref{eq:cnuepvpart} are the portions
    (in green) of the corresponding full circles above the plane
    tangent to the collision surface (the excluded parts of those
    circles are shown in white).}
  \label{fig:collsurface}
\end{figure}

Substituting the parametrisation \eqref{eq:v1v2v3ep} of the velocity
vector $\mathbf{v}\in\bb{S}{2}$ and evaluating its scalar product
with the normal vector \eqref{eq:n}, the velocity integral in equation
\eqref{eq:cnuep} splits into two contributions, integrated over an
arc-length proportional to the angle along the arcs:

\begin{align}
  \int_{2\times \bb{S}{1}: \, \mathbf{v}\cdot\mathbf{n} \geq 0}
  \mathrm{d}\mathbf{v}\,
  &
    (\mathbf{v}\cdot\mathbf{n})
    =
    \frac{1}{\sqrt2}
    \int_{\bb{S}{1}: \, \sqrt{2\ep} \geq \sqrt{1 - 2\ep} \cos\alpha}
    \mathrm{d}\alpha\,
    (\sqrt{2\ep} - \sqrt{1 - 2\ep} \cos\alpha)
    \cr
  &
    \quad
    +
    \frac{1}{\sqrt2}
    \int_{\bb{S}{1}: \,  \sqrt{2\ep} \leq -\sqrt{1 - 2\ep} \cos\alpha}
    \mathrm{d}\alpha\,
    (-\sqrt{2\ep} - \sqrt{1 - 2\ep} \cos\alpha)
    \, ;
  \label{eq:cnuepvpart}
\end{align}
see \fref{fig:collsurface} for a graphical representation. Two separate regimes arise.

On the one hand, when the piston's energy is less than the ball's, the
condition $\sqrt{2\ep} \geq \sqrt{1 - 2\ep} \cos \alpha$ in the first
of the two integrals on the right-hand side of equation
\eqref{eq:cnuepvpart} is equivalent to

\begin{equation}
  \label{eq:alphaabove}
  \arccos \sqrt{\frac{\ep}{\tfrac{1}{2} - \ep}} \leq \alpha
  \leq 2 \pi - \arccos \sqrt{\frac{\ep}{\tfrac{1}{2} - \ep}}
  \,.
\end{equation}

\noindent
Performing the integral, we obtain the contribution,

\begin{align}
  \frac{1}{\sqrt2}
  \int_{\bb{S}{1}: \, \sqrt{2\ep} \geq \sqrt{1 - 2\ep}
  \cos\alpha}
  &
    \mathrm{d}\alpha\,
    (\sqrt{2\ep} - \sqrt{1 - 2\ep} \cos\alpha)
    \cr
  &
    =
    2\sqrt{\ep} \left( \pi - \arccos \sqrt{\frac{\ep}{\tfrac{1}{2} -
    \ep}} \right)
    + 2 \sqrt{ \tfrac{1}{2} - 2 \ep}
    \,.
  \label{eq:intalphaabove}
\end{align}

Likewise, the condition $\sqrt{2\ep} \leq - \sqrt{1 - 2\ep} \cos
\alpha$ in the second of the two integrals on the right-hand side of
equation \eqref{eq:cnuepvpart} is equivalent to

\begin{equation}
  \label{eq:alphabelow}
  \pi - \arccos \sqrt{\frac{\ep}{\tfrac{1}{2} - \ep}} \leq \alpha
  \leq \pi + \arccos \sqrt{\frac{\ep}{\tfrac{1}{2} - \ep}}
  \,.
\end{equation}

\noindent
Thus the second integral yields the contribution

\begin{align}
  \frac{1}{\sqrt2}
  \int_{\bb{S}{1}: \, \sqrt{2\ep} \leq - \sqrt{1 - 2\ep}
  \cos\alpha}
  &
    \mathrm{d}\alpha\,
    (-\sqrt{2\ep} - \sqrt{1 - 2\ep} \cos\alpha)
    \cr
  &
    =
    - 2\sqrt{\ep} \arccos \sqrt{\frac{\ep}{\tfrac{1}{2} - \ep}}
    + 2 \sqrt{ \tfrac{1}{2} - 2 \ep}
    \,.
  \label{eq:intalphabelow}
\end{align}

The contribution to equation \eqref{eq:cnuep} from the velocity
integral in the corresponding energy interval is thus given by the sum
of equations \eqref{eq:intalphaabove} and \eqref{eq:intalphabelow}.


On the other hand, when the piston's energy is larger than the ball's,
the condition $\sqrt{2\ep} \geq \sqrt{1 - 2\ep} \cos \alpha$ in the
first of the two integrals on the right-hand side of equation
\eqref{eq:cnuepvpart} holds true for all angles $\alpha$. The result
of the integration,

\begin{equation}
      \frac{1}{\sqrt2}
    \int_{\bb{S}{1}: \, \sqrt{2\ep} \geq \sqrt{1 - 2\ep} \cos\alpha}
    \mathrm{d}\alpha\,
    (\sqrt{2\ep} - \sqrt{1 - 2\ep} \cos\alpha)
    = 2 \pi \sqrt{\ep}
    \,,
\end{equation}

\noindent yields the only contribution to the velocity integral in
equation \eqref{eq:cnuep}.

\begin{acknowledgements} 
  The authors gratefully acknowledge fruitful discussions with Tam\'as
  Tasn\'ady. They are also grateful to Makiko Sasada for communicating
  her results prior to publishing. TG and PN wish to acknowledge the
  hospitality of the Institute of Mathematics at the Budapest
  University of Technology and Economics where part of this work was
  conducted. TG also wishes to acknowledge stimulating discussions
  with the DinAmicI community during their 2015 workshop, held in
  Corinaldo, Italy.  The authors are grateful to the Erwin
  Schr\"odinger Institute, Vienna, for their hospitality on the
  occasion of the workshop \emph{Mixing Flows and Averaging Methods}
  during which this work was finalized. Finally, they also express
  their gratitude to the referees for valuable comments. 
  PB, DSz and IPT acknowledge
  the financial support of the Hungarian National Foundation for
  Scientific Research (OTKA): grant  T104745 and of Stiftung Aktion
  \"Osterreich-Ungarn: grant 87\"ou6. TG is financially supported
  by the (Belgian) FRS-FNRS.
\end{acknowledgements}


\end{document}